\documentclass{amsart}
\usepackage{amssymb,latexsym}
\usepackage{epsfig}
\usepackage{eufrak}
\usepackage{amsmath}
\usepackage{mathrsfs}

\theoremstyle{plain}
\newtheorem{theorem}{Theorem}

\newtheorem{lemma}{Lemma}

\theoremstyle{definition}

\theoremstyle{remark}

\numberwithin{equation}{section}

\newtheorem{remark}{Remark}

\newcommand{\be}{\begin{equation}}
\newcommand{\ee}{\end{equation}}
\newcommand{\bee}{\begin{eqnarray*}}
\newcommand{\eee}{\end{eqnarray*}}
\newcommand{\bel}{\begin{eqnarray}}
\newcommand{\eel}{\end{eqnarray}}
\newcommand{\bec}{\begin{cases}}
\newcommand{\eec}{\end{cases}}
\newcommand{\bem}{\begin{bmatrix}}
\newcommand{\eem}{\end{bmatrix}}

\newcommand{\bed}{\begin{description}}
\newcommand{\eed}{\end{description}}
\newcommand{\bei}{\begin{itemize}}
\newcommand{\eei}{\end{itemize}}
\newcommand{\ben}{\begin{enumerate}}
\newcommand{\een}{\end{enumerate}}

\newcommand{\beL}{\begin{lemma}}
\newcommand{\eeL}{\end{lemma}}
\newcommand{\beT}{\begin{theorem}}
\newcommand{\eeT}{\end{theorem}}

\newcommand{\bpf}{\begin{pf}}
\newcommand{\epf}{\end{pf}}

\newcommand{\pfbox}{\hfill\mbox{$\Box$}}
\newenvironment{pf}{\paragraph*{Proof{\rm.}}}{\pfbox\bigskip}

\begin{document}

\title[Binomial Confidence Interval]
{Explicit Formula for Constructing Binomial Confidence Interval with Guaranteed Coverage Probability}

\author{Xinjia Chen, Kemin Zhou and Jorge L. Aravena}

\address{Department of Electrical and Computer Engineering\\
Louisiana State University\\
Baton Rouge, LA 70803}

\email{chan@ece.lsu.edu\\
kemin@ece.lsu.edu\\
aravena@ece.lsu.edu}

\thanks{This
research was supported in part by grants from NASA (NCC5-573) and LEQSF (NASA /LEQSF(2001-04)-01).}


\date{June 2006}

\begin{abstract}

In this paper, we derive an explicit formula for constructing the confidence
interval of binomial parameter with guaranteed coverage probability.
The formula overcomes the limitation of
normal approximation which is asymptotic
in nature and thus inevitably introduce unknown errors in applications.
Moreover, the formula is very tight in comparison with classic Clopper-Pearson's
approach  from the perspective of interval width.  Based on the rigorous formula,
we also obtain approximate formulas with excellent performance of coverage probability.

\end{abstract}

\maketitle

\section{Classic Confidence Intervals}

The construction of confidence interval of binomial parameter is frequently encountered
in communications and many other areas of science and engineering.
Clopper and Pearson \cite{Clo} has provided a rigorous approach for constructing confidence interval.
However, the computational complexity involved with this approach is very high.
The standard technique is to use normal approximation which is not accurate for rare events,
especially in the context of studying the bit error rate of communication systems,
blocking probability of communication networks and probability of
instability of uncertain dynamic systems.  Moreover, it has been recently proven by
Brown, Cai and DasGupta \cite{BCD1, BCD2} that the standard normal
approximation approach is persistently poor.  The coverage probability of the confidence interval
can be significantly below the specified confidence level even for very large sample sizes.
Since in many situations, it is desirable to
quickly construct a confidence interval with guaranteed coverage probability,
our goal is to derive a simple and rigorous formula for confidence interval construction.

Let the probability space be denoted as $(\Omega, F,P)$ where
$\Omega, F, P$ are the sample space, the algebra of events and
the probability measure respectively. Let $X$ be a Bernoulli
random variable with distribution ${\rm Pr} \{ X=1\} =
\mathbb{P}_X, \;\;{\rm Pr} \{X=0\} = 1-\mathbb{P}_X$ where
$\mathbb{P}_X \in (0,1)$.   Let the sample size $N$ and confidence
parameter $\delta \in (0,1)$ be fixed. We refer an observation
with value 1 as a successful trial. Let $K$ denote the number of
successful trials during the $N$ i.i.d. sampling experiments. Let
$k = K(\omega)$ where $\omega$ is a sample point in the sample
space $\Omega$.

\subsection{Clopper-Pearson Confidence Limits}

The classic Clopper-Pearson lower confidence limit $L_{N,k,\delta}$ and
   upper confidence limit $U_{N,k,\delta}$ are given respectively by
\[
L_{N,k,\delta} \stackrel{\mathrm{def}}{=} \left\{\begin{array}{ll}
   0 \;\;\;&  {\rm if}\; k=0\\
   \underline{p} \; \;\;\;&
   {\rm if}\; k > 0
\end{array} \right.\;\;\;{\rm and}\;\;\;
U_{N,k,\delta} \stackrel{\mathrm{def}}{=} \left\{\begin{array}{ll}
   1  &  {\rm if}\; k=N\\
   \overline{p} \;
   & {\rm if}\; k < N
\end{array} \right.
\]
where $\underline{p} \in (0,1)$ is the solution of the following equation
\be
\label{clp1}
\sum_{j=0}^{k-1} {N \choose j}
   \underline{p}^j (1- \underline{p})^{N-j} = 1-\frac{\delta}{2}
\ee
and $\overline{p} \in
(0,1)$ is the solution of the following equation
\be
\label{clp2}
\sum_{j=0}^{k} {N \choose j}
   \overline{p}^j (1- \overline{p})^{N-j} = \frac{\delta}{2}.
\ee
The probabilistic implication of the confidence limits can be illustrated as follows:
Define random variable $L: \Omega \rightarrow [0,1]$ by
$L (\omega) = L_{N,K(\omega),\delta} \;\;\forall \omega \in \Omega$
and random variable $U: \Omega \rightarrow [0,1]$ by
$U (\omega) = U_{N,K(\omega),\delta} \;\;\forall \omega \in \Omega$.
Then
\[
{\rm Pr}\{ L \leq \mathbb{P}_X \leq U\} > 1-\delta.
\]
The exact value of ${\rm Pr}\{ L \leq \mathbb{P}_X \leq U\}$ is
referred as the {\it coverage probability}. Accordingly, we refer
${\rm Pr}\{ \mathbb{P}_X < L\; {\rm or} \; \mathbb{P}_X > U\}$ as
the {\it error probability}.

\subsection{Normal Approximation}

It is easy to see that the equations~(\ref{clp1}) and~(\ref{clp2}) are very
hard to solve and thus the confidence limits are very difficult
 to determine using Clopper-Pearson's approach.
For large sample size, it is computationally
prohibitive.  To get around the difficulty, normal
approximation has been widely used to develop simple approximate formulas
(see, for example, \cite{BCD1, BCD2, Hald, JKW} and the references therein ).
The basis of the normal approximation is the Central Limit Theorem, i.e.,
\[
\lim_{N \rightarrow \infty} \Pr \left\{  \frac{ \left|
\frac{K}{N} - \mathbb{P}_X \right|} { \sqrt{ \frac{ \mathbb{P}_X
(1-\mathbb{P}_X) } {N} }   }  < z \right\}
 = 2 \Phi(z) - 1
\]
where $z > 0$ and $\Phi(.)$ is the normal distribution function.
Let $Z_{\frac{\delta}{2}}$ be the critical value
such that $\Phi( Z_{\frac{\delta}{2}} ) = 1 - \frac{\delta}{2}$.
It follows that
\[
\lim_{N \rightarrow \infty} \Pr \left\{   \frac{K}{N}   -
Z_{\frac{\delta}{2}} \sqrt{ \frac{ \mathbb{P}_X (1-\mathbb{P}_X) }
{N} } < \mathbb{P}_X <
 \frac{K}{N}  + Z_{\frac{\delta}{2}} \sqrt{ \frac{ \mathbb{P}_X (1-\mathbb{P}_X) } {N} } \right\}  = 1 - \delta,
\]
i.e.,

\vspace{0.1in}

$\lim_{N \rightarrow \infty} \;
\Pr \left\{  \frac{ \frac{K}{N}  + \frac{Z_{\frac{\delta}{2}}^2}{2N}
 - Z_{\frac{\delta}{2}} \sqrt{ \frac{ \frac{K}{N} (1-\frac{K}{N}) } {N} +  \frac{Z_{\frac{\delta}{2}}^2}{4N^2} } }
{1 + \frac{Z_{\frac{\delta}{2}}^2}{N} } < \mathbb{P}_X < \frac{
\frac{K}{N}  + \frac{Z_{\frac{\delta}{2}}^2}{2N}
 + Z_{\frac{\delta}{2}} \sqrt{ \frac{ \frac{K}{N} (1-\frac{K}{N}) } {N} +  \frac{Z_{\frac{\delta}{2}}^2}{4N^2} } }
{ 1 + \frac{Z_{\frac{\delta}{2}}^2}{N}   } \right\}  = 1 - \delta$.

\vspace{0.1in}

Since $\frac{Z_{\frac{\delta}{2}}^2}{N}  \approx 0$ for
sufficiently large sample size $N$, the lower and upper
confidence limits can be estimated respectively as
\[
\widetilde{L} \approx  \frac{k}{N}  - Z_{\frac{\delta}{2}} \sqrt{
\frac{ \frac{k}{N} (1-\frac{k}{N}) } {N} } \]
 and
 \[
 \widetilde{U}
\approx \frac{k}{N} + Z_{\frac{\delta}{2}} \sqrt{ \frac{
\frac{k}{N} (1-\frac{k}{N}) } {N} }.
\]

The critical problem with the normal approximation is that it is of
 asymptotic nature.  It is not clear
 how large the sample size is sufficient for the approximation error to be negligible.
Such an asymptotic approach is not good enough for many practical applications involving rare events.

\section{Rigorous Formula}

It is desirable to have a simple formula
which is rigorous and very tight for the confidence interval construction.
We now propose the following simple formula for constructing the confidence limits.

\begin{theorem} \label{Tape_Massart}
Define \be \mathcal{ L} (k) \stackrel{\mathrm{def}}{=}
\frac{k}{N} + \frac{3}{4} \; \frac{ 1 - \frac{2k}{N} - \sqrt{ 1 +
4 \theta \; k ( 1- \frac{k}{N}) } } {1 + \theta N }, \quad k=0,1,
\cdots, N \label{CI_l} \ee and \be \mathcal{ U} (k)
\stackrel{\mathrm{def}}{=} \frac{k}{N} + \frac{3}{4} \; \frac{ 1
- \frac{2k}{N} + \sqrt{ 1 + 4 \theta \; k ( 1- \frac{k}{N}) } }
{1 + \theta N }, \quad k=0,1, \cdots, N \label{CI_u} \ee with
$\theta = \frac{9}{ 8 \ln \frac{2}{\delta} }$. Then $\Pr
\left\{   \mathcal{ L} (K)< \mathbb{P}_X < \mathcal{ U} (K)
\right\}
> 1 - \delta$. Moreover,
\[
\mathcal{ L}(k) < L_{N, k, \delta} < U_{N, k, \delta} < \mathcal{
U}(k).
\]

\end{theorem}

\begin{remark}
$\mathcal{ L}(k)$ and $\mathcal{ U}(k)$ are tight bounds
for the classic Clopper-Pearson confidence limits
$L_{N, k, \delta}$ and $U_{N, k, \delta}$ (See Figures 1-12).  A bisection search can be performed
based on such bounds for computing the classic Clopper-Pearson confidence limits.
\end{remark}


To show Theorem~\ref{Tape_Massart}, we need some preliminary results.
The following Lemma \ref{Massart} is due to Massart \cite{Massart}.

\begin{lemma} \label{Massart}
${\Pr} \left\{ \frac{K}{N} \geq \mathbb{P}_X + \epsilon \right\}
\leq \exp \left(  - \frac{N \epsilon^2} {2 ( \mathbb{P}_X +
\frac{\epsilon}{3} )\;( 1 - \mathbb{P}_X - \frac{\epsilon}{3} )}
\right)$ for all  $\epsilon \in (0, 1 - \mathbb{P}_X)$.
\end{lemma}

Of course, the above upper bound holds trivially for $\epsilon
\geq 1 - \mathbb{P}_X$. Thus, Lemma~\ref{Massart} is actually
true for any $\epsilon >0$.

\begin{lemma} \label{Massart2} ${\Pr} \left\{ \frac{K}{N} \leq \mathbb{P}_X - \epsilon \right\} \leq
\exp \left(  - \frac{N \epsilon^2} {2 ( \mathbb{P}_X -
\frac{\epsilon}{3} )\;( 1 - \mathbb{P}_X + \frac{\epsilon}{3} )}
\right)$ for all  $\epsilon > 0$.
\end{lemma}

\begin{pf}

Define $Y = 1 - X$.  Then $\mathbb{P}_Y = 1 - \mathbb{P}_X$.  At
the same time when we are conducting $N$ i.i.d. experiments for
$X$, we are also conducting $N$ i.i.d. experiments for $Y$. Let
the number of successful trials of the experiments for $Y$ be
denoted as $K_Y$.  Obviously, $K_Y = N - K$. Applying
Lemma~\ref{Massart} to $Y$, we have
\[
{\Pr} \left\{ \frac{K_Y}{N} \geq \mathbb{P}_Y + \epsilon \right\}
\leq \exp \left(  - \frac{N \epsilon^2} {2 ( \mathbb{P}_Y +
\frac{\epsilon}{3} )\;( 1 - \mathbb{P}_Y - \frac{\epsilon}{3} )}
\right).
\]
It follows that
\[
{\Pr} \left\{ \frac{N - K}{N} \geq 1 - \mathbb{P}_X + \epsilon
\right\} \leq \exp \left(  - \frac{N \epsilon^2} {2 ( 1 -
\mathbb{P}_X + \frac{\epsilon}{3} )\;[ 1 - (1- \mathbb{P}_X) -
\frac{\epsilon}{3} ]} \right).
\]
The proof is thus completed by observing that ${\Pr} \left\{
\frac{N - K}{N} \geq 1 - \mathbb{P}_X + \epsilon \right\} = {\Pr}
\left\{ \frac{K}{N} \leq \mathbb{P}_X - \epsilon \right\}$.
\end{pf}

The following lemma can be found in \cite{david}.
\begin{lemma} \label{decrease} $\sum_{j=0}^{k} {N \choose j}
   x^j (1- x)^{N-j}$ decreases monotonically with respect to $x \in (0,1)$ for
$k=0,1,\cdots,N$.
\end{lemma}

\begin{lemma} \label{lem_massart} $ \sum_{j=0}^{k} {N \choose j}
   x^j (1- x)^{N-j} \leq \exp \left( - \frac{ N ( x - \frac{k}{N} )^2 }
{ 2 \; ( \frac{2}{3} x + \frac{k}{3N} ) \; ( 1- \frac{2}{3} x - \frac{k}{3N} ) } \right)
\;\;\forall x \in (\frac{k}{N}, 1)$ for
$k=0,1,\cdots,N$.
\end{lemma}
{\bf Proof.}  Consider binomial random variable $X$ with
parameter $\mathbb{P}_X
> \frac{k}{N}$.  Let $K$ be the number
 of successful trials
 during $N$ i.i.d. sampling experiments.  Then
\[
\sum_{j=0}^{k} {N \choose j}
   \mathbb{P}_X^j (1- \mathbb{P}_X)^{N-j} = {\rm Pr} \{ K \leq k\}.
   \]
   Note that ${\rm Pr} \{ K \leq k\} = {\rm Pr} \left\{ \frac{K}{N} \leq \mathbb{P}_X - \left( \mathbb{P}_X -
\frac{k}{N} \right) \right\}$. Applying Lemma \ref{Massart2} with
$\epsilon = \mathbb{P}_X - \frac{k}{N} > 0$, we have
\begin{eqnarray*}
\sum_{j=0}^{k} {N \choose j}
   \mathbb{P}_X^j (1- \mathbb{P}_X)^{N-j}
& \leq & \exp \left(  - \frac{N ( \mathbb{P}_X - \frac{k}{N} )^2}
{2 ( \mathbb{P}_X - \frac{\mathbb{P}_X - \frac{k}{N}}{3} )\;( 1 -
\mathbb{P}_X + \frac{\mathbb{P}_X - \frac{k}{N}}{3} )}
\right)\\
& = & \exp \left( - \frac{ N ( \mathbb{P}_X  - \frac{k}{N} )^2 }
{ 2 \; ( \frac{2}{3} \mathbb{P}_X  + \frac{k}{3N} ) \; ( 1-
\frac{2}{3} \mathbb{P}_X  - \frac{k}{3N} ) } \right).
\end{eqnarray*}
Since the argument holds for arbitrary binomial random variable
$X$ with $\mathbb{P}_X
> \frac{k}{N}$, the proof of the lemma is thus completed.
$\;\;\;\;\square$

\begin{lemma}\label{lem_Massart_B}
$\sum_{j=0}^{k-1} {N \choose j}
   x^j (1- x)^{N-j}
   \geq 1 - \exp \left( - \frac{ N ( x - \frac{k}{N} )^2 }
{ 2 \; ( \frac{2}{3} x + \frac{k}{3N} ) \; ( 1- \frac{2}{3} x - \frac{k}{3N} ) } \right)
\;\;
\forall x \in (0,\frac{k}{N})$ for
$k=1,\cdots,N$.
\end{lemma}
{\bf Proof.}  Consider binomial random variable $X$ with parameter
$ \mathbb{P}_X < \frac{k}{N}$.  Let $K$ be the number
 of successful trials
 during $N$ i.i.d. sampling experiments.
 Then
 \[
\sum_{j=0}^{k-1} {N \choose j}
   \mathbb{P}_X^j (1- \mathbb{P}_X)^{N-j}
 =  {\rm Pr} \{K < k\} =
 {\rm Pr} \left\{ \frac{K}{N} < \mathbb{P}_X  +
 (\frac{k}{N} - \mathbb{P}_X) \right\}.
 \]
 Applying Lemma \ref{Massart} with
$\epsilon = \frac{k}{N} - \mathbb{P}_X > 0$, we have that
\begin{eqnarray*}
\sum_{j=0}^{k-1} {N \choose j}
   \mathbb{P}_X^j (1- \mathbb{P}_X)^{N-j}
& \geq & 1 - \exp \left(  - \frac{N (\frac{k}{N} - \mathbb{P}_X
)^2} {2 ( \mathbb{P}_X + \frac{ \frac{k}{N} - \mathbb{P}_X }{3}
)\; ( 1 - \mathbb{P}_X - \frac{ \frac{k}{N} - \mathbb{P}_X }{3} )}
\right)\\
& = & 1 - \exp \left( - \frac{ N ( \mathbb{P}_X  - \frac{k}{N}
)^2 } { 2 \; ( \frac{2}{3} \mathbb{P}_X  + \frac{k}{3N} ) \; ( 1-
\frac{2}{3} \mathbb{P}_X  - \frac{k}{3N} ) } \right).
\end{eqnarray*}
Since the argument holds for arbitrary binomial random variable
$X$ with $\mathbb{P}_X < \frac{k}{N}$, the proof of the lemma is
thus completed. $\;\;\;\;\square$

\begin{lemma}\label{lem8} Let $0 \leq k \leq N$. Then
$L_{N,k,\delta} < U_{N,k,\delta}$.
\end{lemma}
{\bf Proof.}  Obviously, the lemma is true for k$=0,N$.
We consider the case
that $1 \leq k \leq N-1$.  Let $\mathcal{ S} (N,k,x) = \sum_{j=0}^{k} {N \choose j}
x^j (1- x)^{N-j}$ for $x \in (0,1)$.  Notice that
\[
\mathcal{ S} (N,k,\overline{p})=
\mathcal{ S} (N,k-1,\overline{p})+{N \choose
k}
   \overline{p}^k (1-\overline{p})^{N-k}=
   \frac{\delta}{2}.
\]
Thus
\[
\mathcal{ S} (N,k-1,\underline{p})-\mathcal{ S} (N,k-1,\overline{p})=
   1-\frac{\delta}{2}-\left[\frac{\delta}{2}-{N \choose k}
   \overline{p}^k (1-\overline{p})^{N-k}\right].
   \]
Notice that $\delta \in (0,1)$ and that $\overline{p} \in (0,1)$,
we have that
\[
\mathcal{ S} (N,k-1,\underline{p})-\mathcal{ S} (N,k-1,\overline{p})
= 1-\delta+
   {N \choose k} \overline{p}^k (1-\overline{p})^{N-k} > 0.
\]
By Lemma~\ref{decrease}, $\mathcal{ S} (N,k-1,x)$
decreases monotonically with respect to $x$,
we have $\underline{p} <
\overline{p}$ and complete the proof of the lemma.   $\;\;\;\;\square$

\bigskip

We are now in the position to prove Theorem~\ref{Tape_Massart}.
It can be easily verified that $U_{N,k,\delta} \leq \mathcal{ U} (k)$ for $k = 0, \; N$.
We need to show that $U_{N,k,\delta} \leq \mathcal{ U} (k)$ for $0 < k < N$.
Straightforward computation shows that $\mathcal{ U} (k)$ is the only root of
equation
\[
\exp \left( - \frac{ N ( x - \frac{k}{N} )^2 }
{ 2 \; ( \frac{2}{3} x + \frac{k}{3N} ) \; ( 1- \frac{2}{3} x - \frac{k}{3N} ) } \right)
= \frac{\delta}{2}
\]
with
respect to $x \in (\frac{k}{N},\infty)$.  There are two cases: $\mathcal{ U} (k) \geq 1$ and $\mathcal{ U} (k) < 1$.
If $\mathcal{ U} (k) \geq 1$ then $U_{N,k,\delta} \leq \mathcal{ U} (k)$
is trivially true.  We only need to consider the case that $\frac{k}{N} < \mathcal{ U} (k) < 1$.
In this case, it follows from Lemma~\ref{lem_massart}
that
\[
\sum_{j=0}^{k} {N \choose j}
   {[\mathcal{ U} (k)]}^j (1- \mathcal{ U} (k))^{N-j} \leq \exp \left( - \frac{ N ( \mathcal{ U} (k) - \frac{k}{N} )^2 }
{ 2 \; ( \frac{2}{3} \mathcal{ U} (k) + \frac{k}{3N} ) \; ( 1- \frac{2}{3} \mathcal{ U} (k) - \frac{k}{3N} ) } \right)
= \frac{\delta}{2}.
\]
Recall that
 \[
\sum_{j=0}^{k} {N \choose j}
   U_{N,k,\delta}^j (1- U_{N,k,\delta})^{N-j} = \frac{\delta}{2},
   \]
 we have
 \[
 \sum_{j=0}^{k} {N \choose j}
   U_{N,k,\delta}^j (1- U_{N,k,\delta})^{N-j} \geq \sum_{j=0}^{k} {N \choose j}
   {[\mathcal{ U} (k)]}^j (1- \mathcal{ U} (k))^{N-j}.
   \]
 Therefore, by Lemma~\ref{decrease}, we have that $U_{N,k,\delta} \leq \mathcal{ U} (k)$ for $0 < k < N$.
Thus, we have shown that $U_{N,k,\delta} \leq q$ for all $k$.

Similarly, by Lemma~\ref{lem_Massart_B} and Lemma~\ref{decrease}, we can show that
 $L_{N,k,\delta} \geq  \mathcal{ L}(k)$.  By Lemma~\ref{lem8}, we have
 $\mathcal{ L}(k) < L_{N, k, \delta} < U_{N, k, \delta} < \mathcal{ U}(k)$.
 Finally, the proof of Theorem~\ref{Tape_Massart} is
 completed by invoking the probabilistic
implication of the Clopper-Pearson confidence interval.

\section{Numerical Experiments and Empirical Formulas}

In comparison with the Clopper-Pearson's approach,
our approach is very tight from the perspective of interval width
(see, for example, Figures 1-12).
Moreover, there is no comparison on
the computational complexity.
Our formula is simple enough for hand calculation.

Our numerical results are in agreement with the discovery made by
Brown, Cai and DasGupta \cite{BCD1, BCD2}. It can be seen from
Figures 21-27 that the coverage probability of confidence
intervals obtained by the standard normal approximation can be
substantially lower than the specified confidence level $1-
\delta$ (This is true even when the condition for applying the
rule of thumb, i.e., $N\mathbb{P}_X (1- \mathbb{P}_X) > 5$, is
satisfied). Moreover, the situation is worse for smaller
confidence parameter $\delta$. See, for example, Figures 25-27,
if one wishes to make an inference with an error frequency less
than one out of $1000$, using the normal approximation can lead
to a frequency of error higher than $100$ out of $1000$. In light
of the excessively high error rate of inference caused by the
normal approximation,  the rigorous formula may be a better
choice. The rigorous formula guarantees the error probability
below the specify level $\delta$. It should be noted that the
rigorous formula is conservative (with actual error probability
around $10\%$ to $20\%$ of the requirement).

It should be noted that by tuning the parameter $\theta$ in the
rigorous formula, one can obtained simple formulas which meet the
specified confidence levels. For example, to construct confidence
interval with confidence parameter $\delta = 0.05, \; 0.01, \;
0.001$, we can simply compute $\mathcal{ L}(k)$ and $\mathcal{ U}
(k)$ defined in Theorem~\ref{Tape_Massart} with $\theta =
\frac{1}{2}, \; \frac{1}{3}, \; \frac{1}{5}$ respectively (The
values of $\theta$ presented here are not optimal. Better
coverage performance can be achieved by a fine tuning of
$\theta$). More specifically,

\vspace{0.1in}

$\Pr \left\{   \frac{K}{N} + \frac{3}{4} \; \frac{ 1 -
\frac{2K}{N} - \sqrt{ 1 + 2 \; K ( 1- \frac{K}{N}) } } {1 +
\frac{N}{2} } < \mathbb{P}_X < \frac{K}{N} + \frac{3}{4} \;
\frac{ 1 - \frac{2K}{N} + \sqrt{ 1 + 2 \; K ( 1- \frac{K}{N}) } }
{1 + \frac{N}{2} } \right\} \approx 0.95; $

\vspace{0.1in}

$ \Pr \left\{   \frac{K}{N} + \frac{3}{4} \; \frac{ 1 -
\frac{2K}{N} - \sqrt{ 1 + \frac{4 \; K }{3}( 1- \frac{K}{N}) } }
{1 + \frac{N}{3} } < \mathbb{P}_X < \frac{K}{N} + \frac{3}{4} \;
\frac{ 1 - \frac{2K}{N} + \sqrt{ 1 + \frac{4\; K}{3} ( 1-
\frac{K}{N}) } } {1 + \frac{N}{3} } \right\} \approx 0.99; $

\vspace{0.1in}

$ \Pr \left\{   \frac{K}{N} + \frac{3}{4} \; \frac{ 1 -
\frac{2K}{N} - \sqrt{ 1 + \frac{4 \; K }{5}( 1- \frac{K}{N}) } }
{1 + \frac{N}{5} } < \mathbb{P}_X < \frac{K}{N} + \frac{3}{4} \;
\frac{ 1 - \frac{2K}{N} + \sqrt{ 1 + \frac{4\; K}{5} ( 1-
\frac{K}{N}) } } {1 + \frac{N}{5} } \right\} \approx 0.999. $

\vspace{0.1in}

Confidence limits computed by these formulas for different $N$
and $\delta$ are depicted by Figures 13-20. It is interesting to
note that, in most situations, the confidence limits computed by
our empirical formulas almost coincide with the corresponding
limits derived by Clopper-Pearson method. The numerical
investigation of the coverage probability of different confidence
intervals is shown in Figures 21-27.  It can be seen that the
empirical formulas have excellent coverage performance.

\begin{figure}[htbp]
\centerline{\psfig{figure=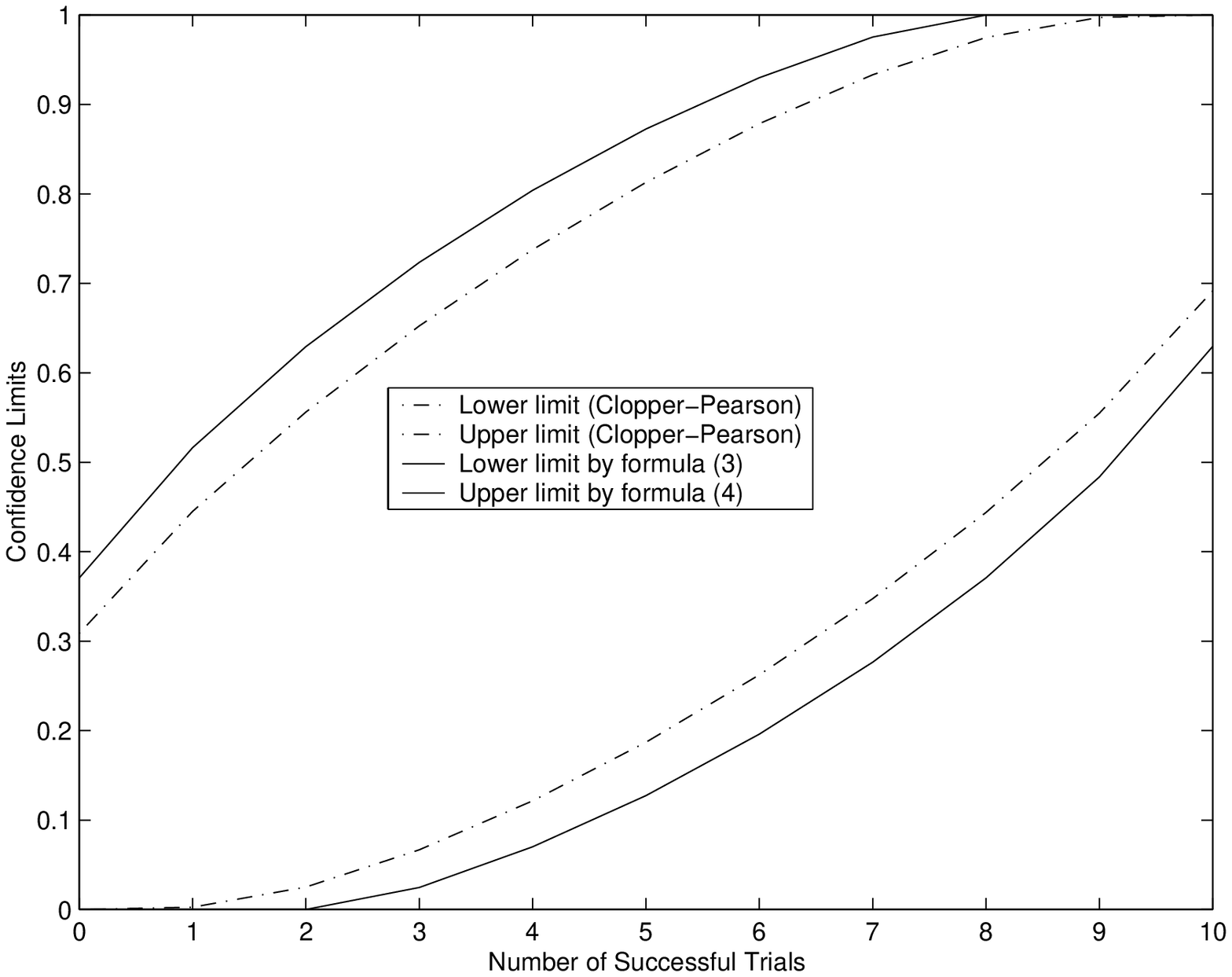, height=3.8in,
width=6.0in }} \caption{ Confidence Interval ($N = 10, \; \delta
= 0.05$.) } \label{fig_1}

\bigskip

\bigskip

\centerline{\psfig{figure=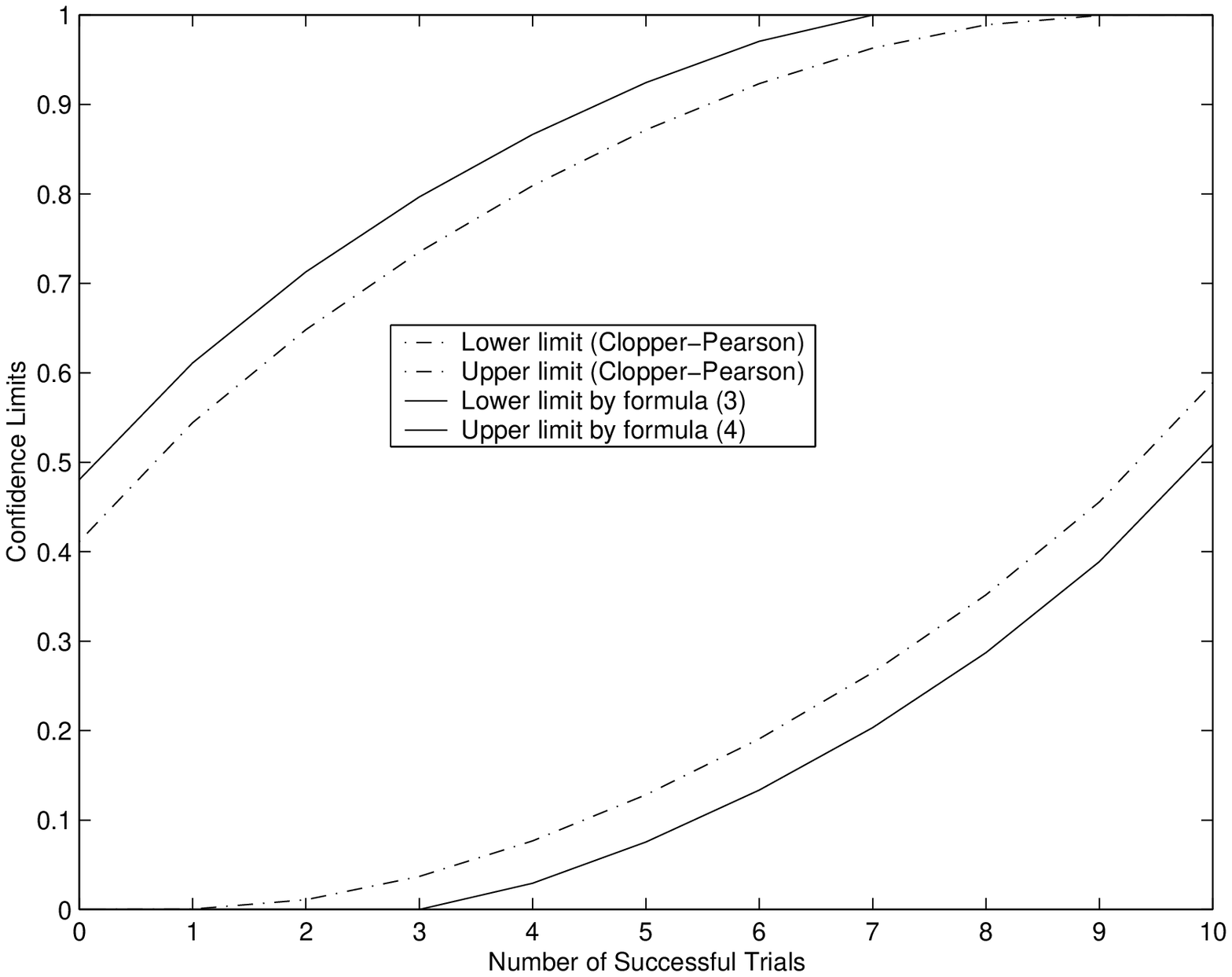, height=3.8in,
width=6.0in }} \caption{ Confidence Interval ($N = 10, \; \delta
= 0.01$.) } \label{fig_2}
\end{figure}

\begin{figure}[htbp]
\centerline{\psfig{figure=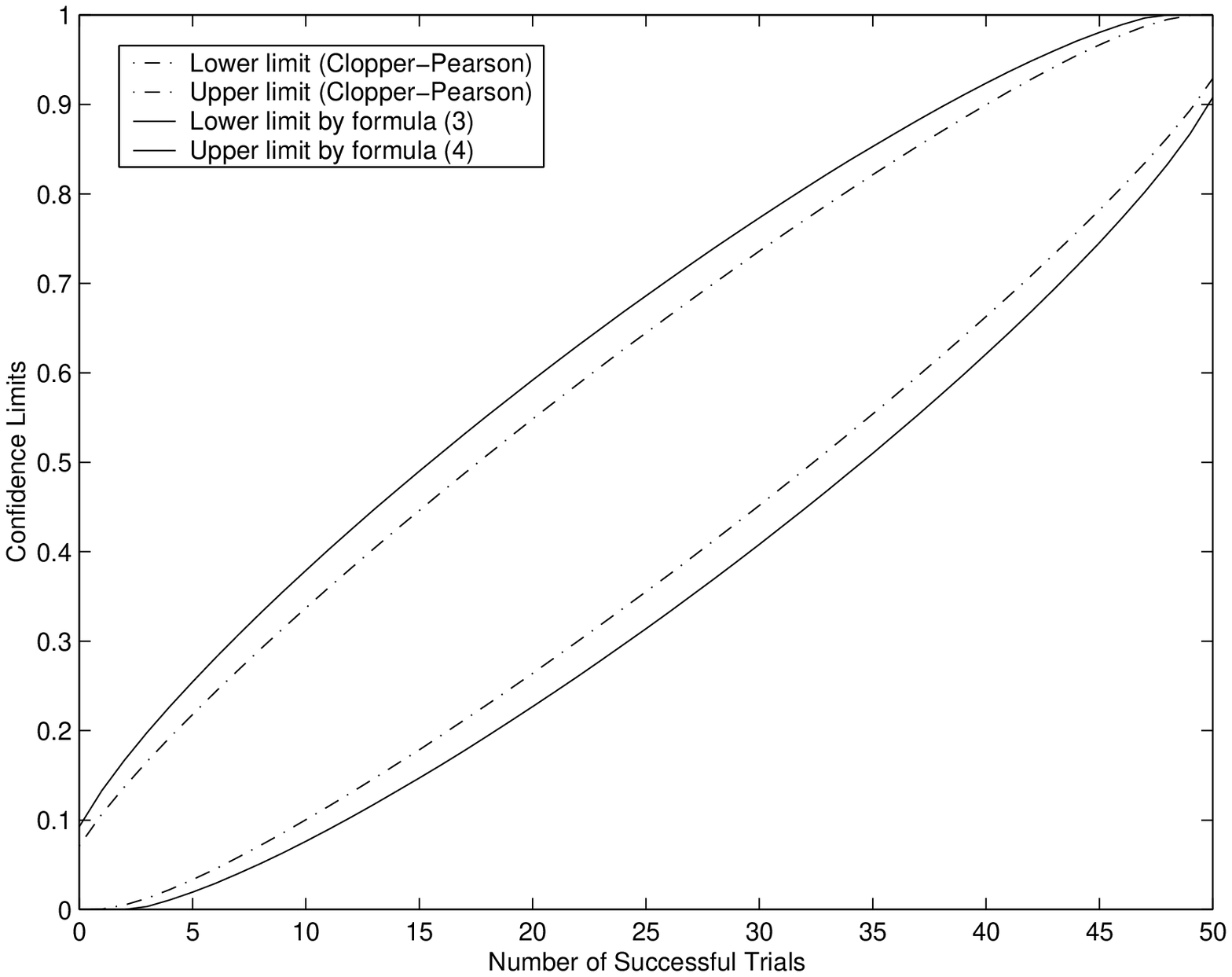, height=3.8in,
width=6.0in }} \caption{ Confidence Interval ($N = 50, \; \delta
= 0.05$.) } \label{fig_3}

\bigskip

\bigskip

\centerline{\psfig{figure=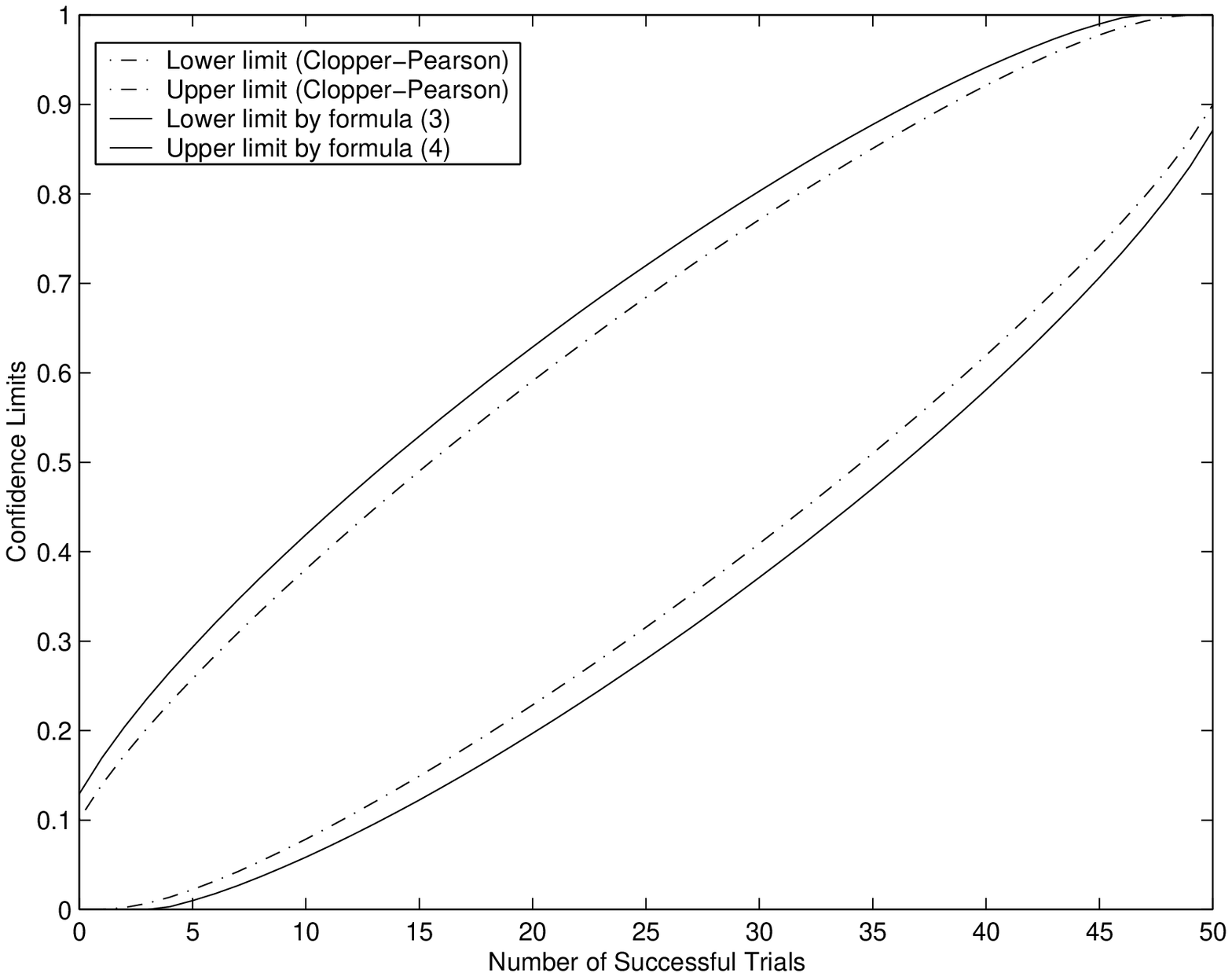, height=3.8in,
width=6.0in }} \caption{ Confidence Interval ($N = 50, \; \delta
= 0.01$.) } \label{fig_4}
\end{figure}

\begin{figure}[htbp]
\centerline{\psfig{figure=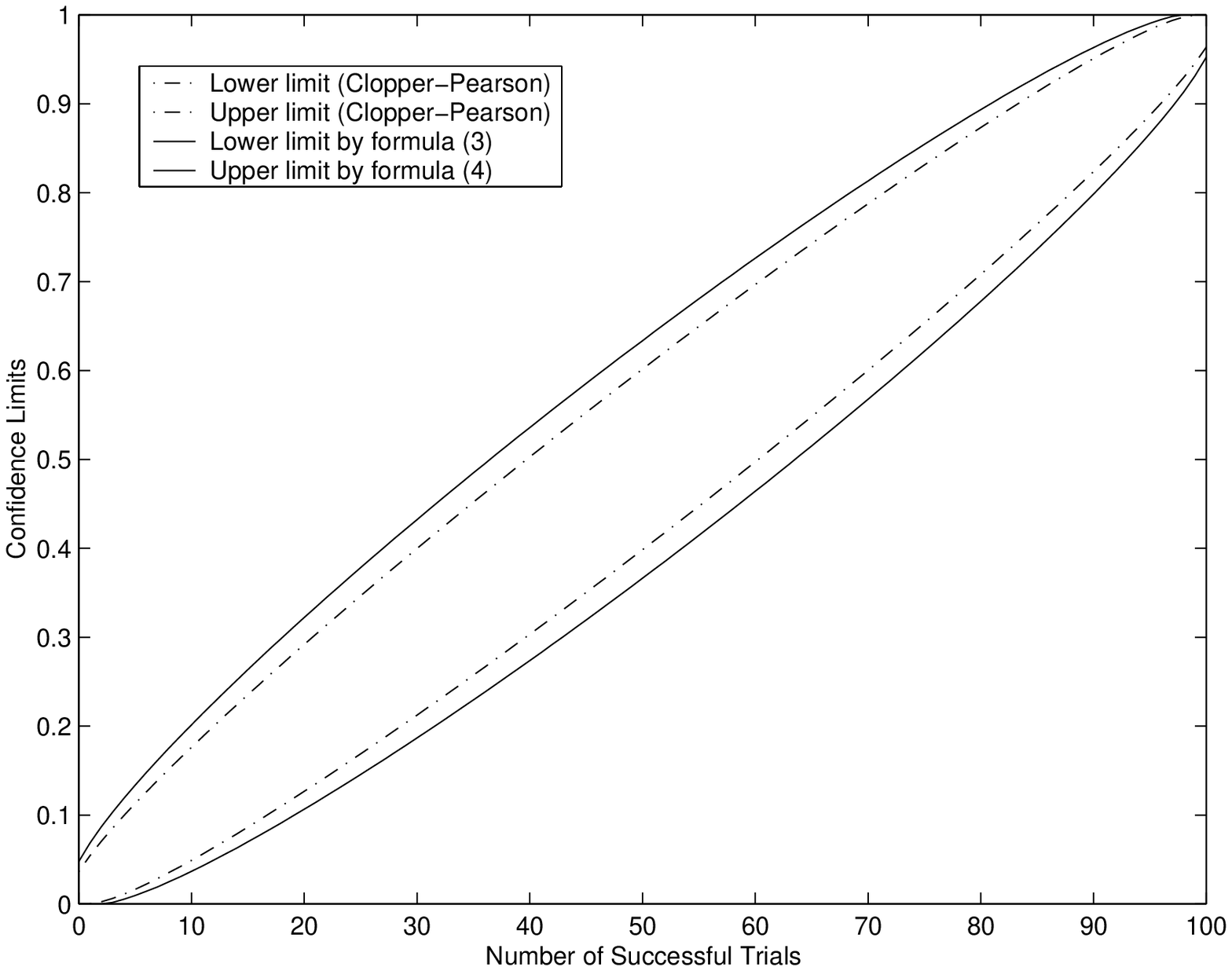, height=3.8in,
width=6.0in }} \caption{ Confidence Interval ($N = 100, \; \delta
= 0.05$.) } \label{fig_5}

\bigskip

\bigskip

\centerline{\psfig{figure=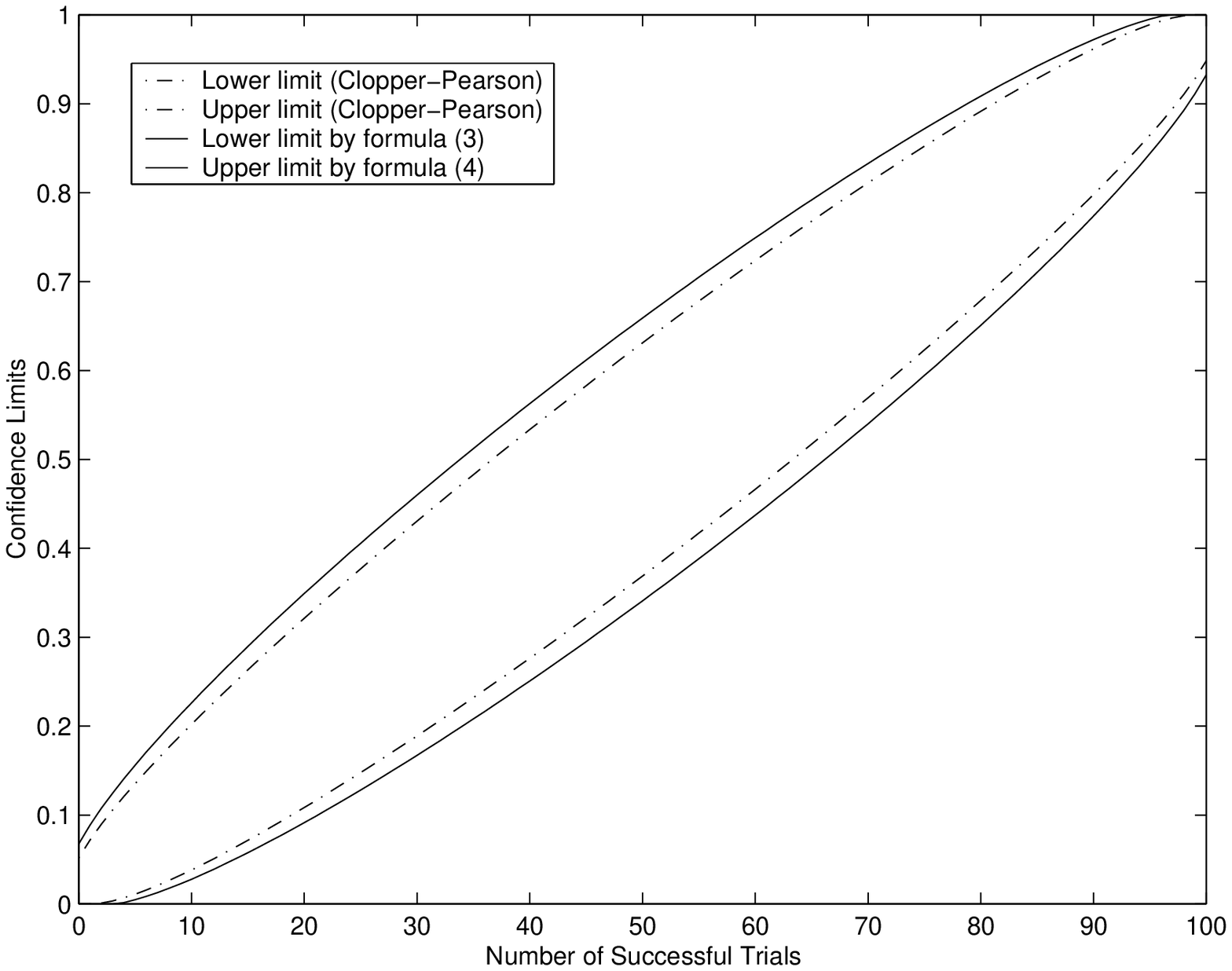, height=3.8in,
width=6.0in }} \caption{ Confidence Interval ($N = 100, \; \delta
= 0.01$.) } \label{fig_6}
\end{figure}

\begin{figure}[htbp]
\centerline{\psfig{figure=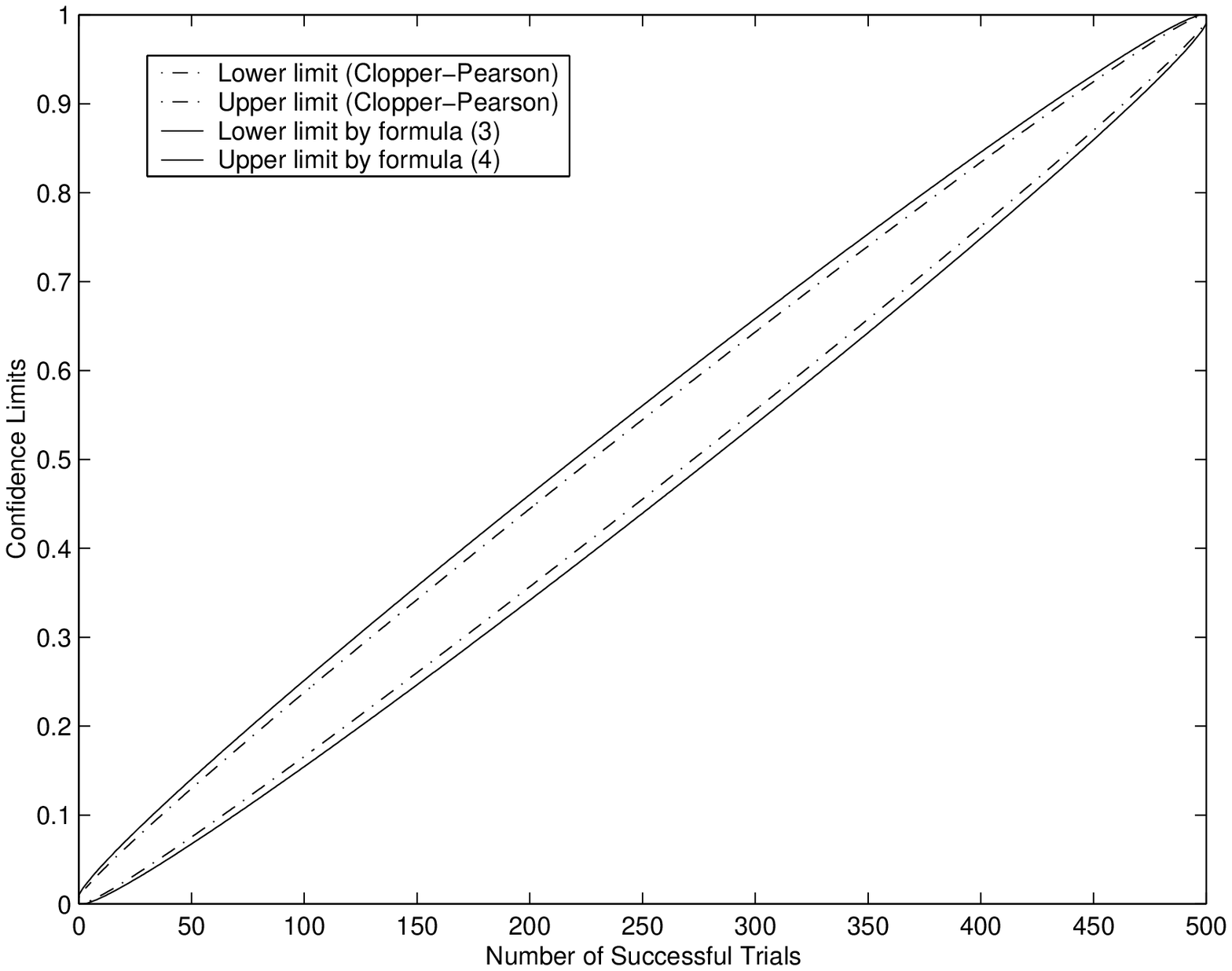, height=3.8in,
width=6.0in }} \caption{ Confidence Interval ($N = 500, \; \delta
= 0.05$.) } \label{fig_7}

\bigskip

\bigskip

\centerline{\psfig{figure=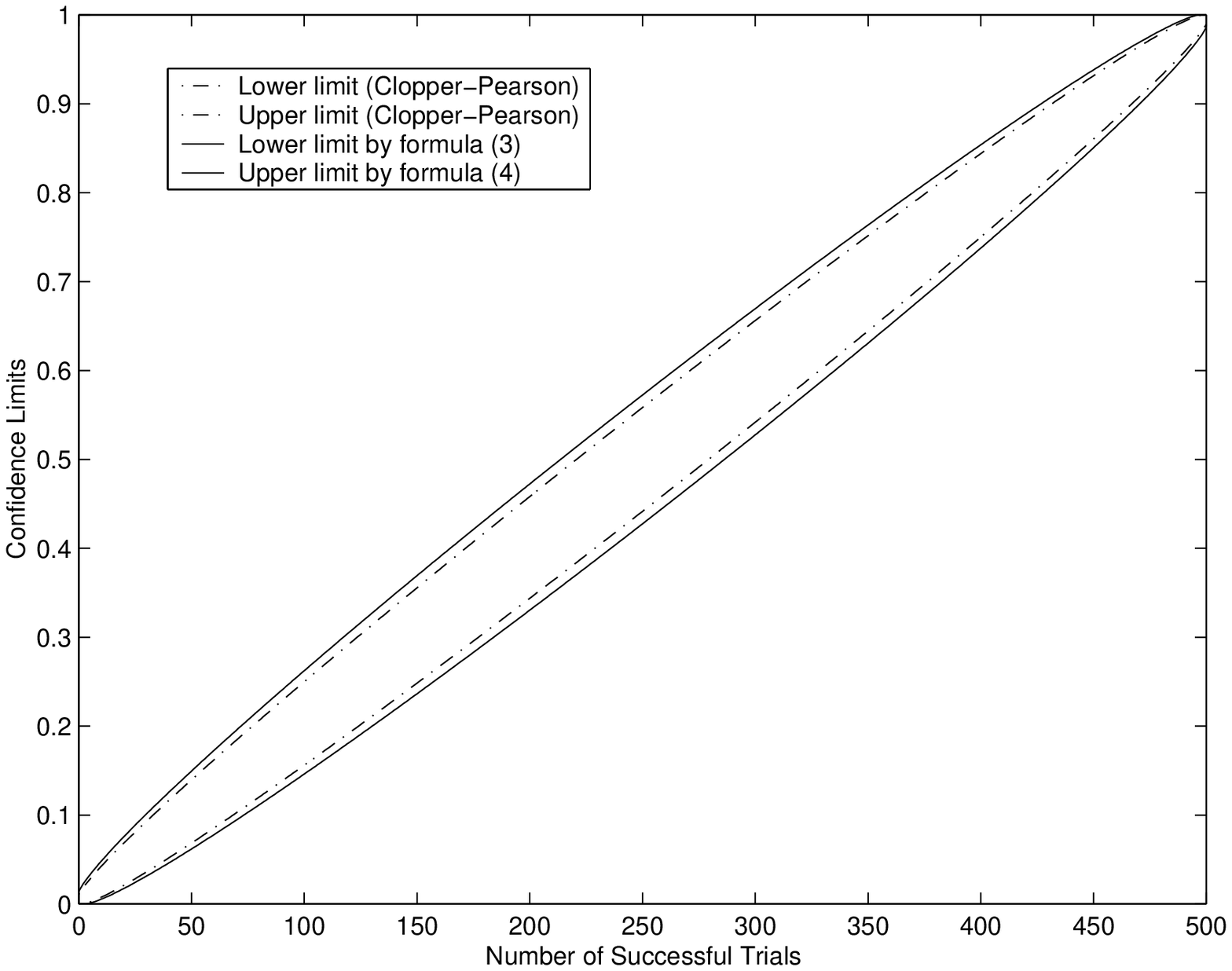, height=3.8in,
width=6.0in }} \caption{ Confidence Interval ($N = 500, \; \delta
= 0.01$.) } \label{fig_8}
\end{figure}

\begin{figure}[htbp]
\centerline{\psfig{figure=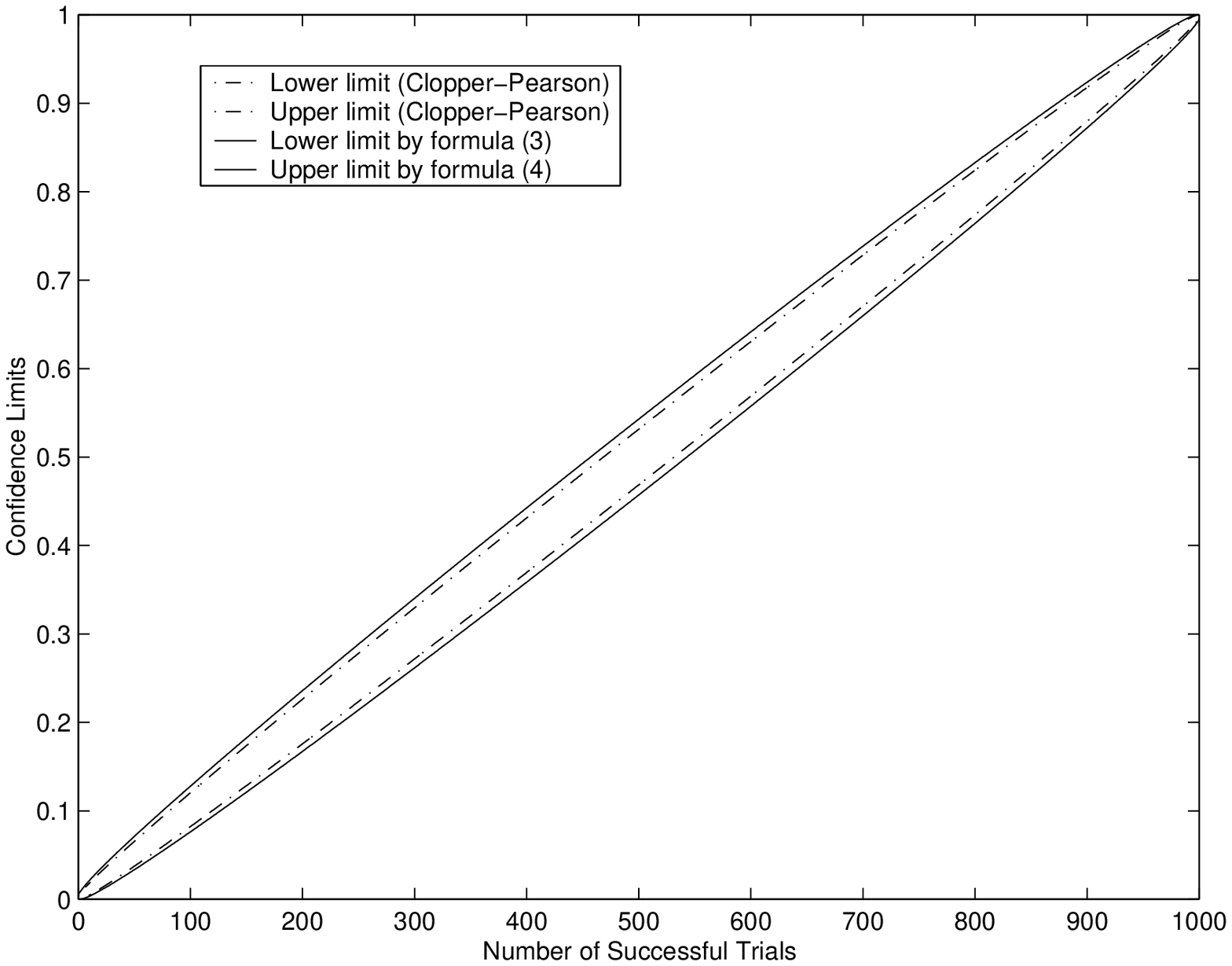, height=3.8in,
width=6.0in }} \caption{ Confidence Interval ($N = 1000, \;
\delta = 0.05$.) } \label{fig_9}

\bigskip

\bigskip

\centerline{\psfig{figure=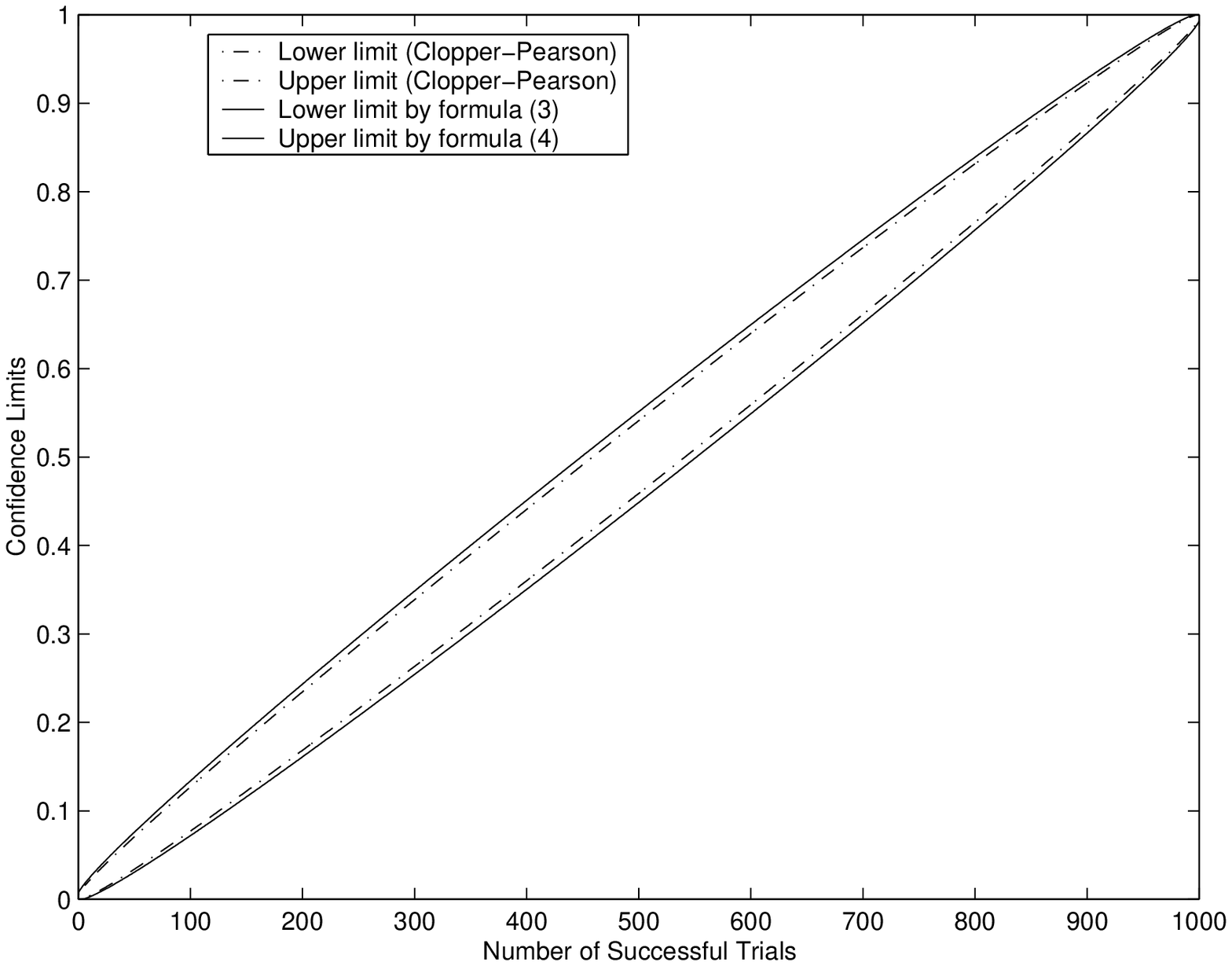, height=3.8in,
width=6.0in }} \caption{ Confidence Interval ($N = 1000, \;
\delta = 0.01$.) } \label{fig_10}
\end{figure}

\begin{figure}[htbp]
\centerline{\psfig{figure=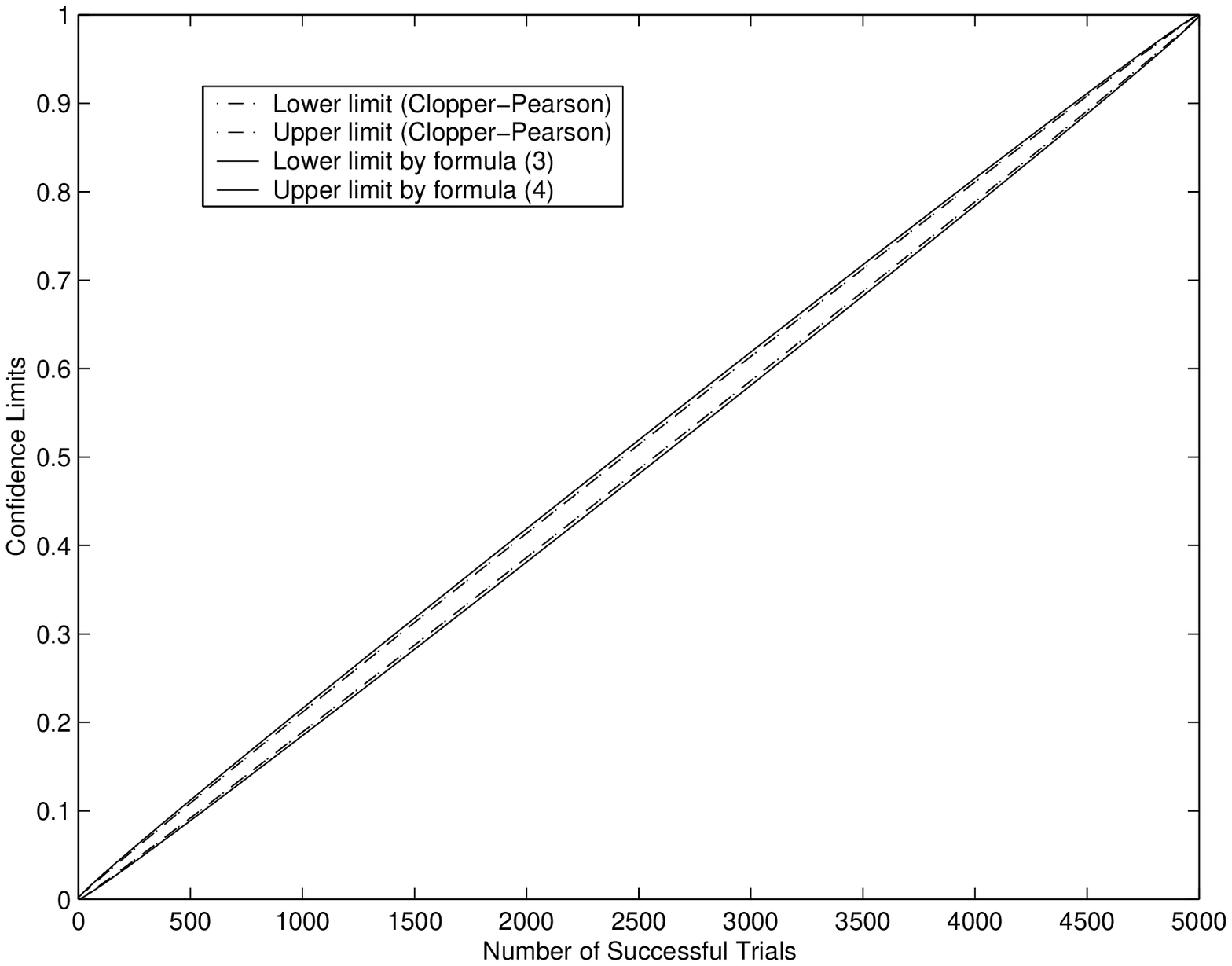, height=3.8in,
width=6.0in }} \caption{ Confidence Interval ($N = 5000, \;
\delta = 0.05$.) } \label{fig_11}

\bigskip

\bigskip

\centerline{\psfig{figure=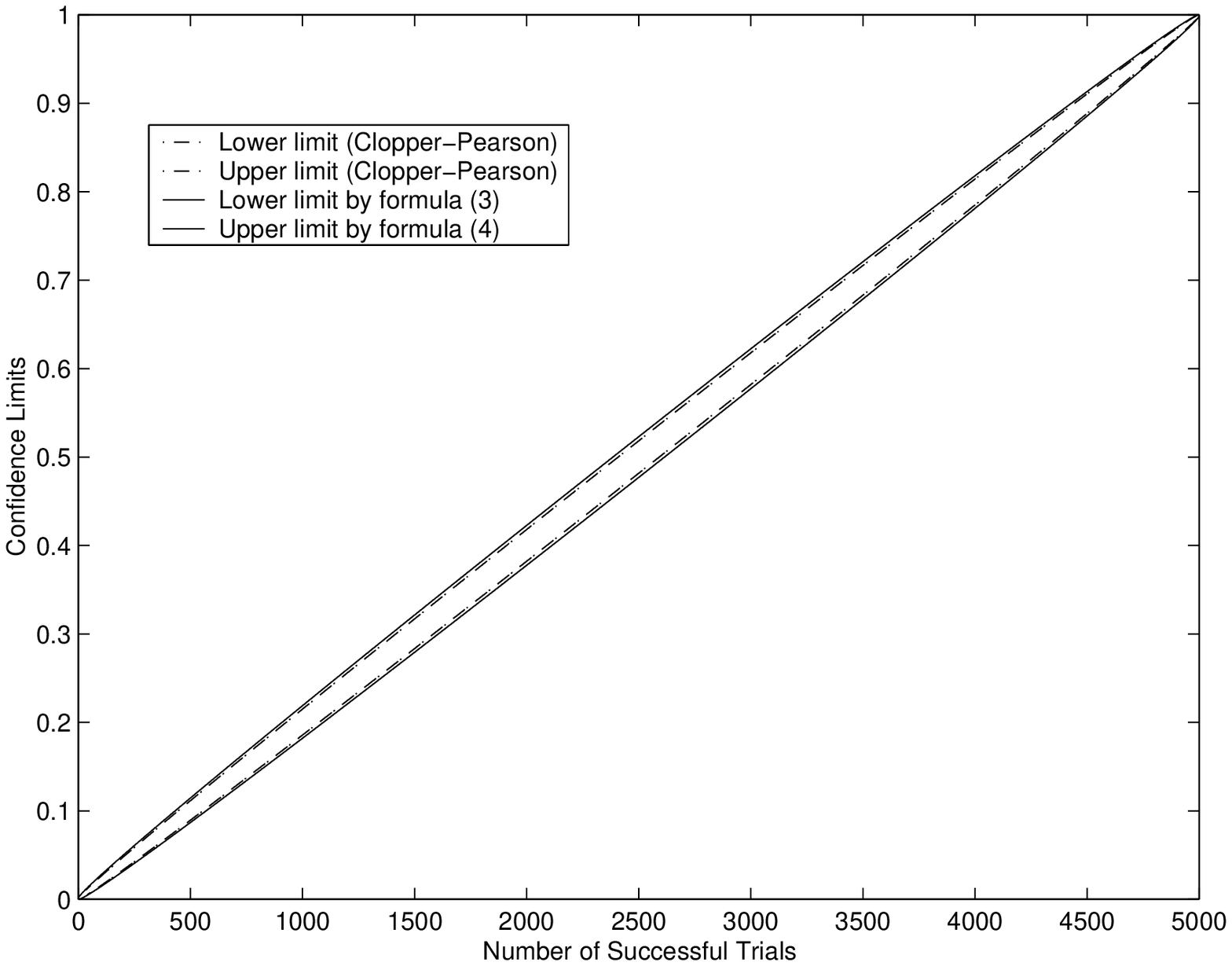, height=3.8in,
width=6.0in }} \caption{ Confidence Interval ($N = 5000, \;
\delta = 0.01$.) } \label{fig_12}
\end{figure}

\begin{figure}[htbp]
\centerline{\psfig{figure=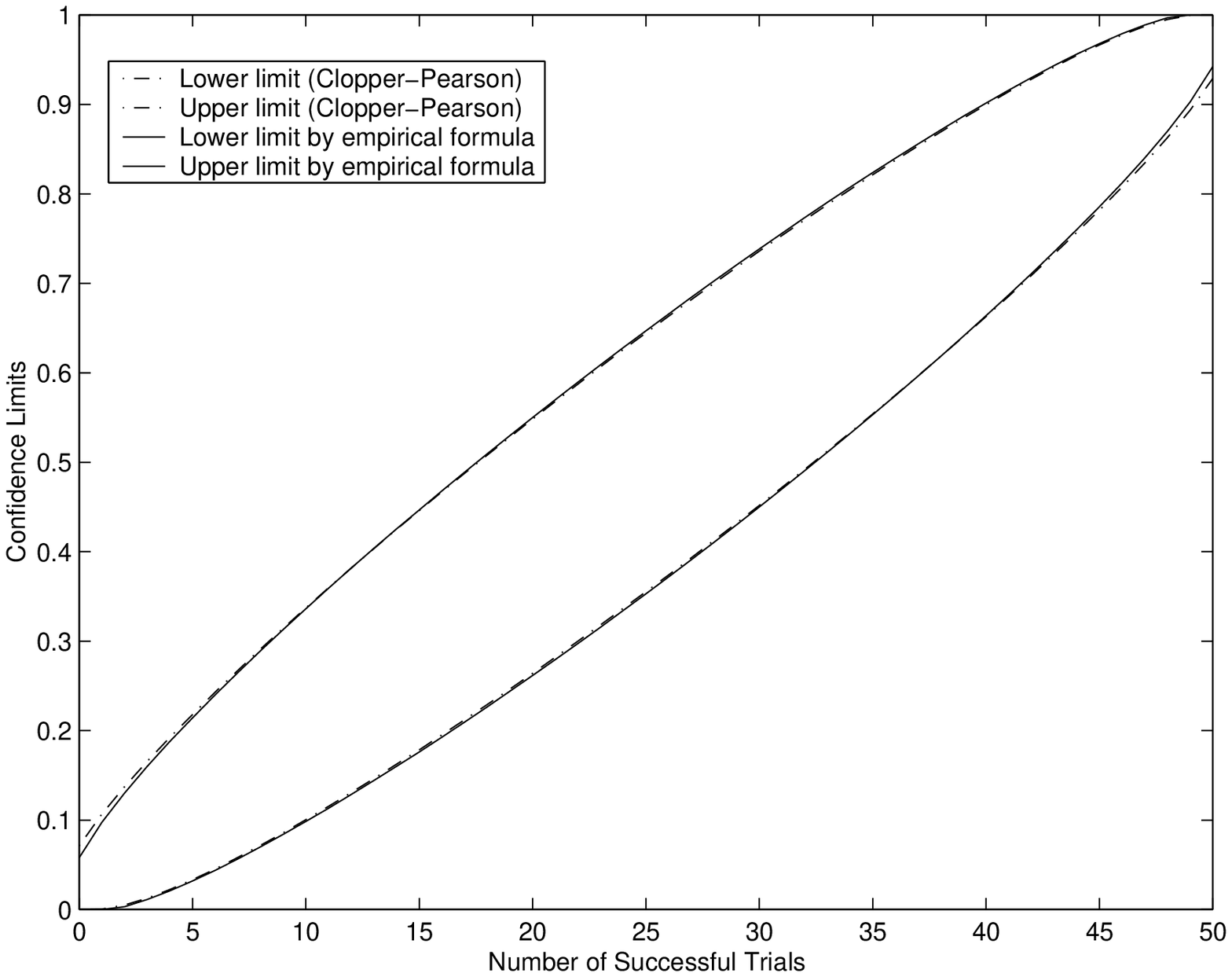, height=3.8in,
width=6.0in }} \caption{ Confidence Interval ($N = 50, \; \delta
= 0.05$.) } \label{fig_13}

\bigskip

\bigskip

\centerline{\psfig{figure=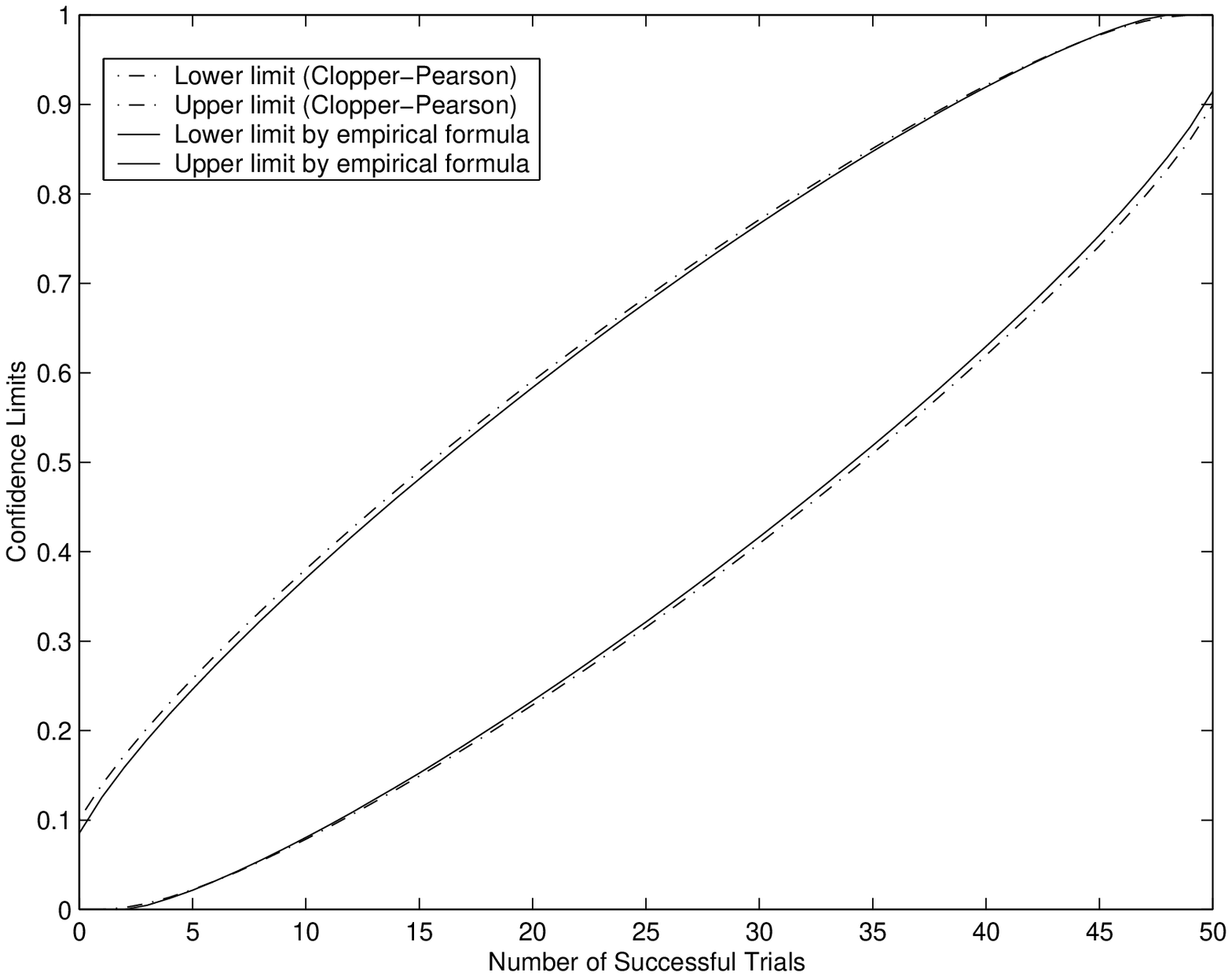, height=3.8in,
width=6.0in }} \caption{ Confidence Interval ($N = 50, \; \delta
= 0.01$.) } \label{fig_14}
\end{figure}

\begin{figure}[htbp]
\centerline{\psfig{figure=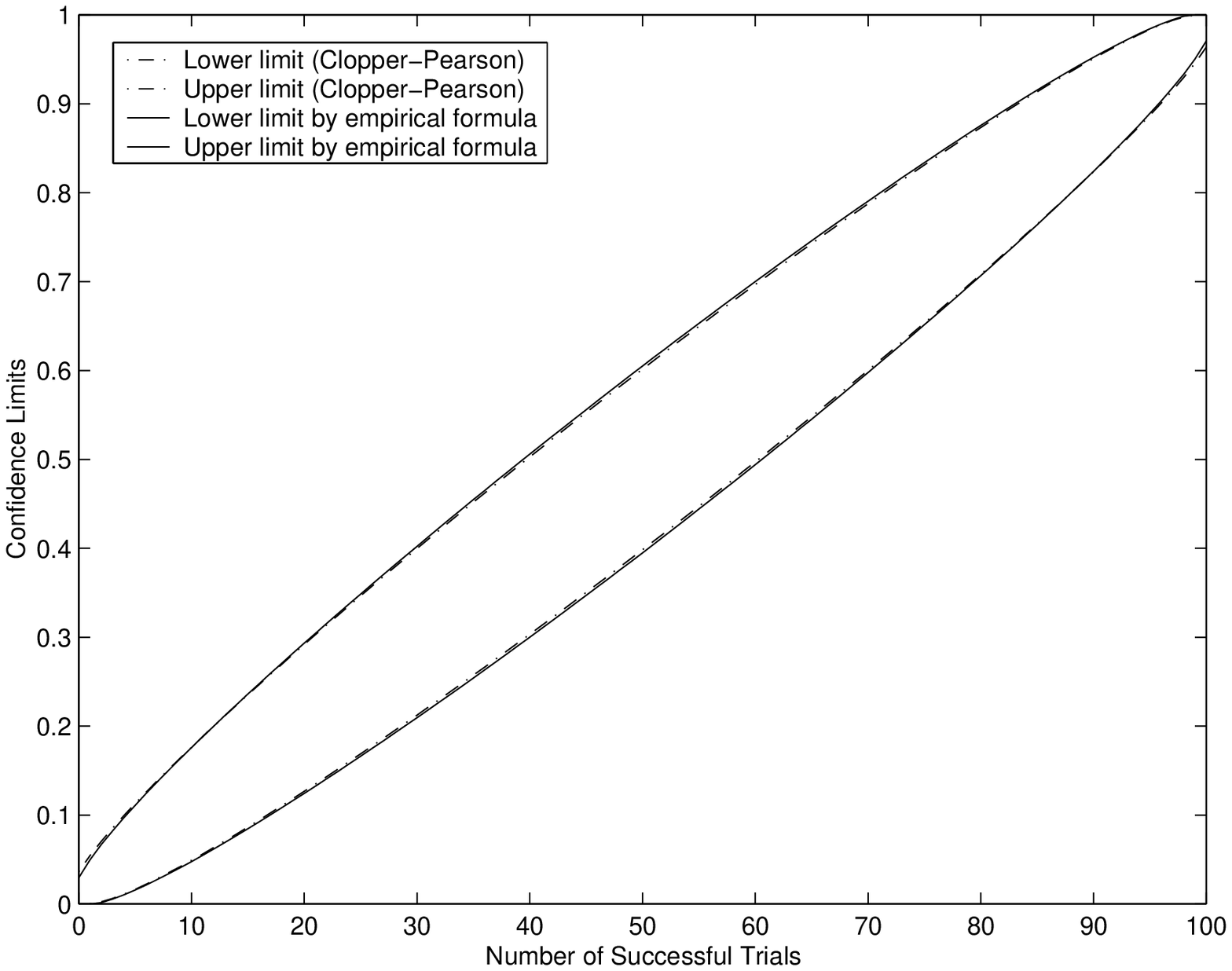, height=3.8in,
width=6.0in }} \caption{ Confidence Interval ($N = 100, \; \delta
= 0.05$.) } \label{fig_15}

\bigskip

\bigskip

\centerline{\psfig{figure=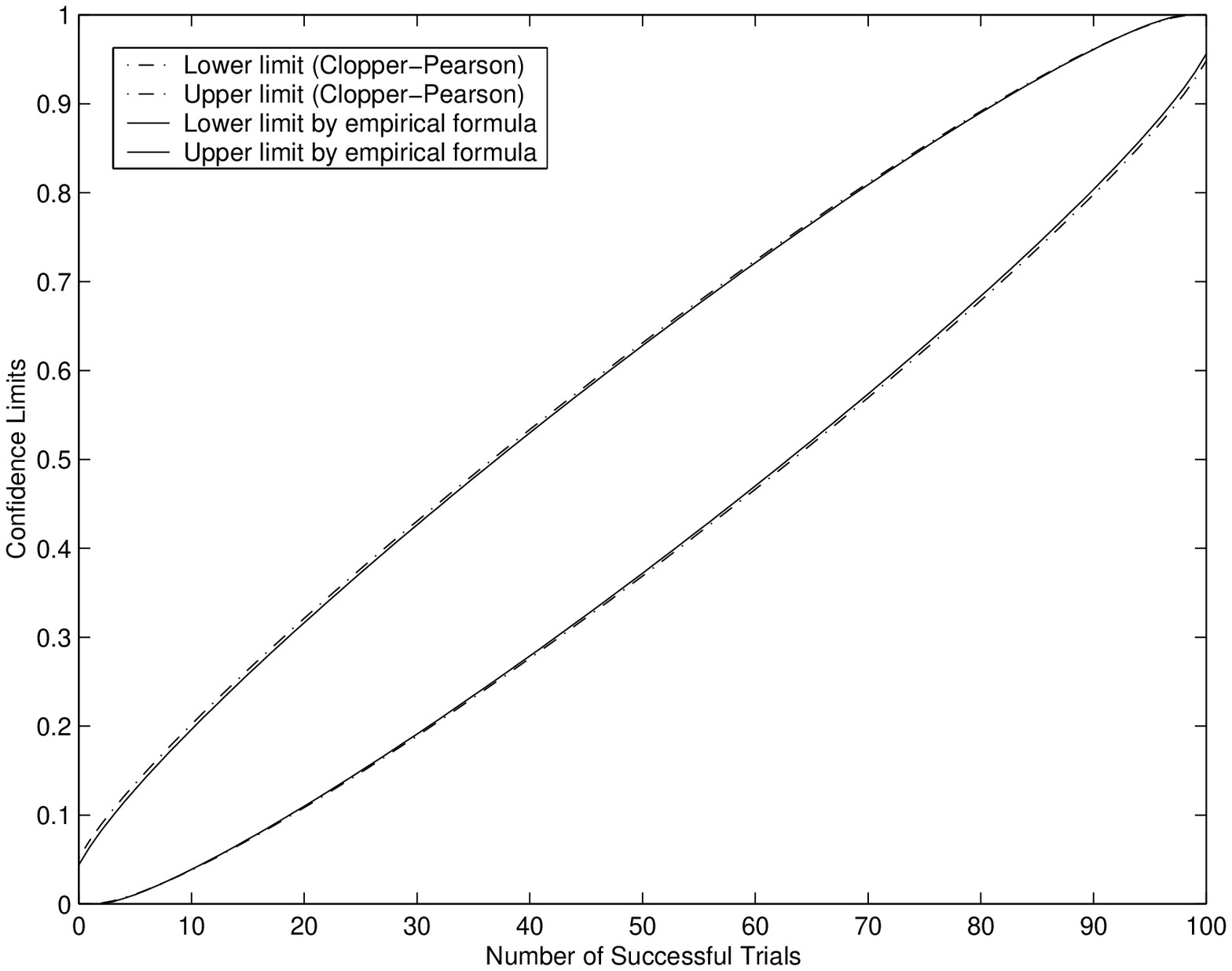, height=3.8in,
width=6.0in }} \caption{ Confidence Interval ($N = 100, \; \delta
= 0.01$.) } \label{fig_16}
\end{figure}

\begin{figure}[htbp]
\centerline{\psfig{figure=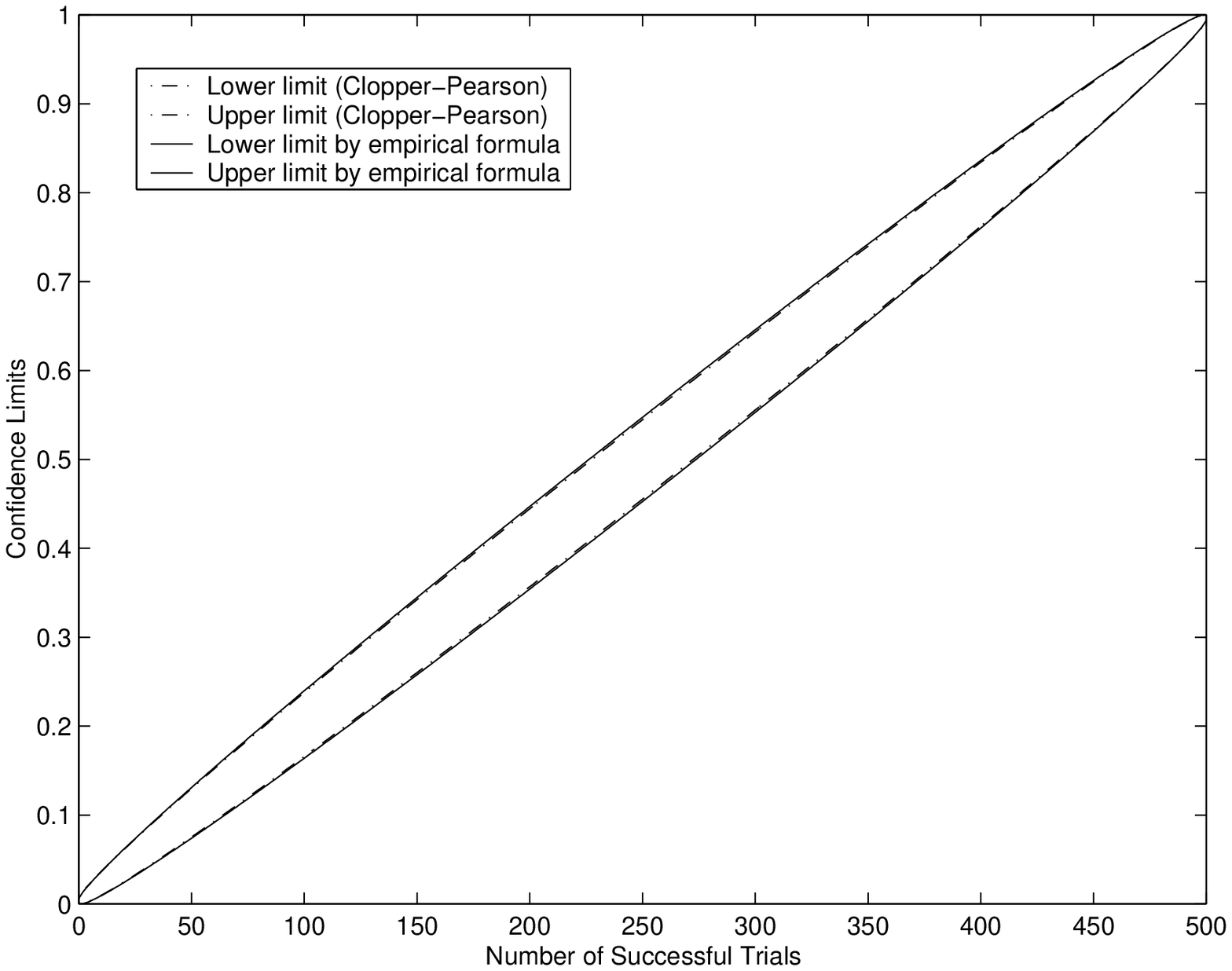, height=3.8in,
width=6.0in }} \caption{ Confidence Interval ($N = 500, \; \delta
= 0.05$.) } \label{fig_17}

\bigskip

\bigskip

\centerline{\psfig{figure=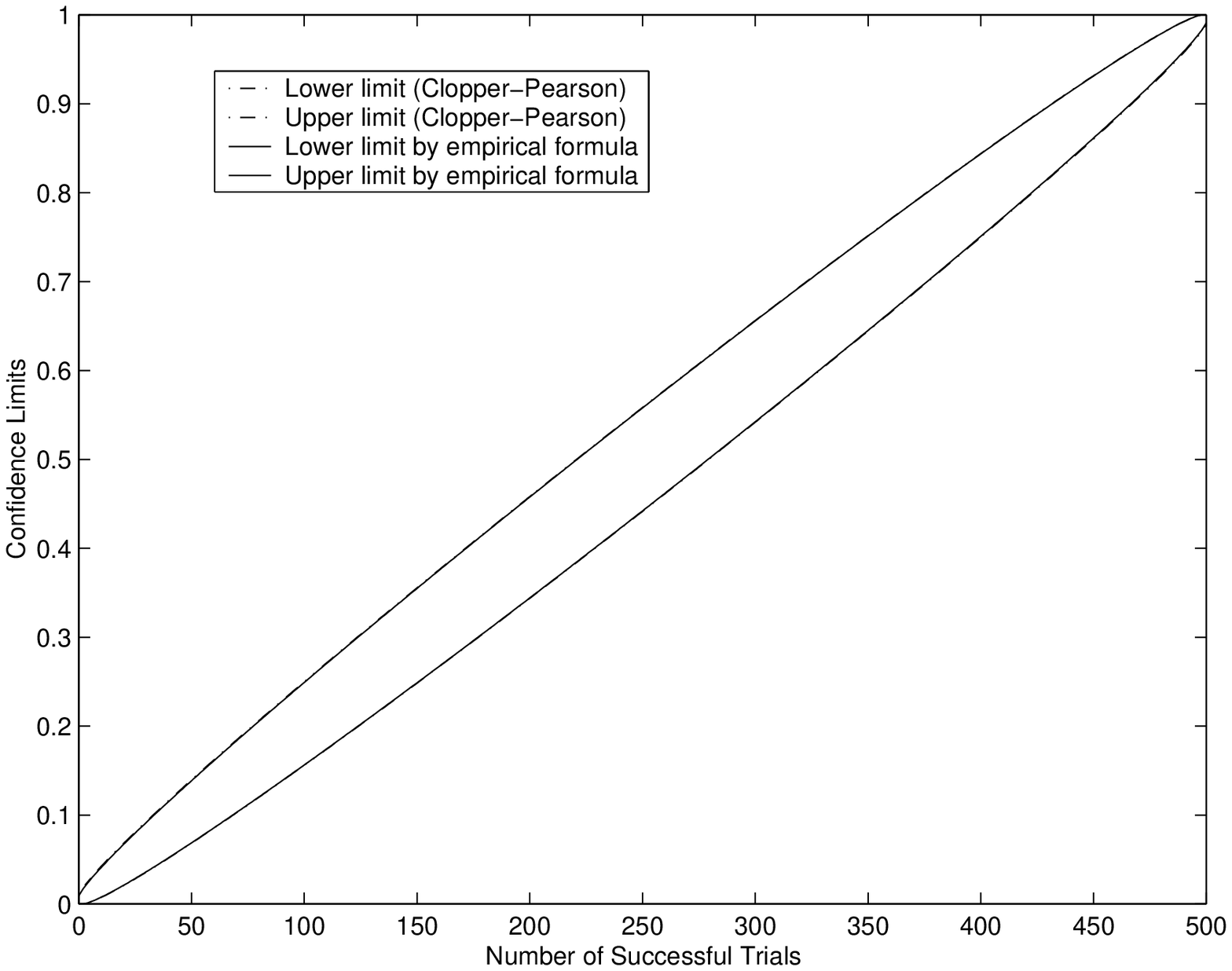, height=3.8in,
width=6.0in }} \caption{ Confidence Interval ($N = 500, \; \delta
= 0.01$.) } \label{fig_18}
\end{figure}

\begin{figure}[htbp]
\centerline{\psfig{figure=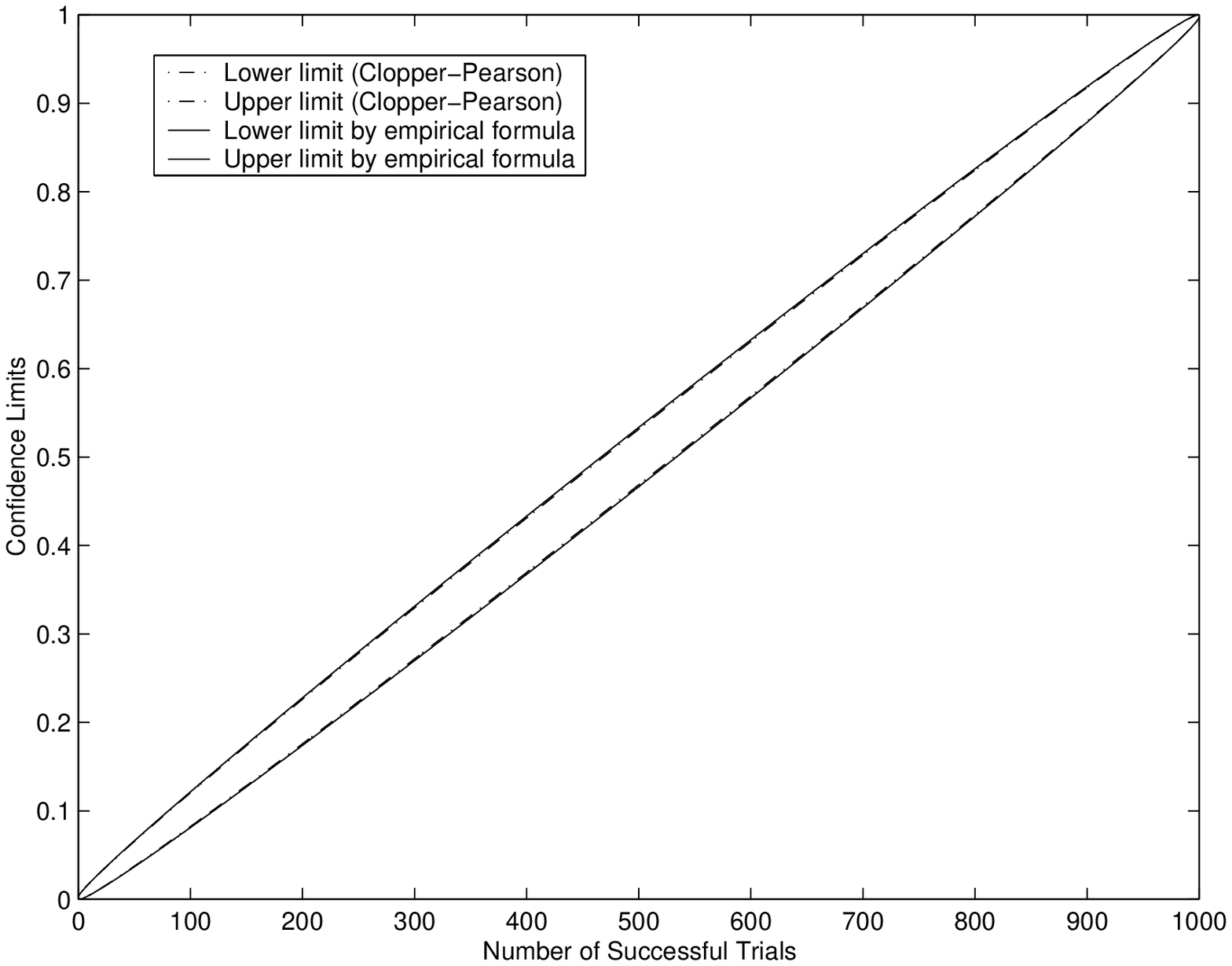, height=3.8in,
width=6.0in }} \caption{ Confidence Interval ($N = 1000, \;
\delta = 0.05$.) } \label{fig_19}

\bigskip

\bigskip

\centerline{\psfig{figure=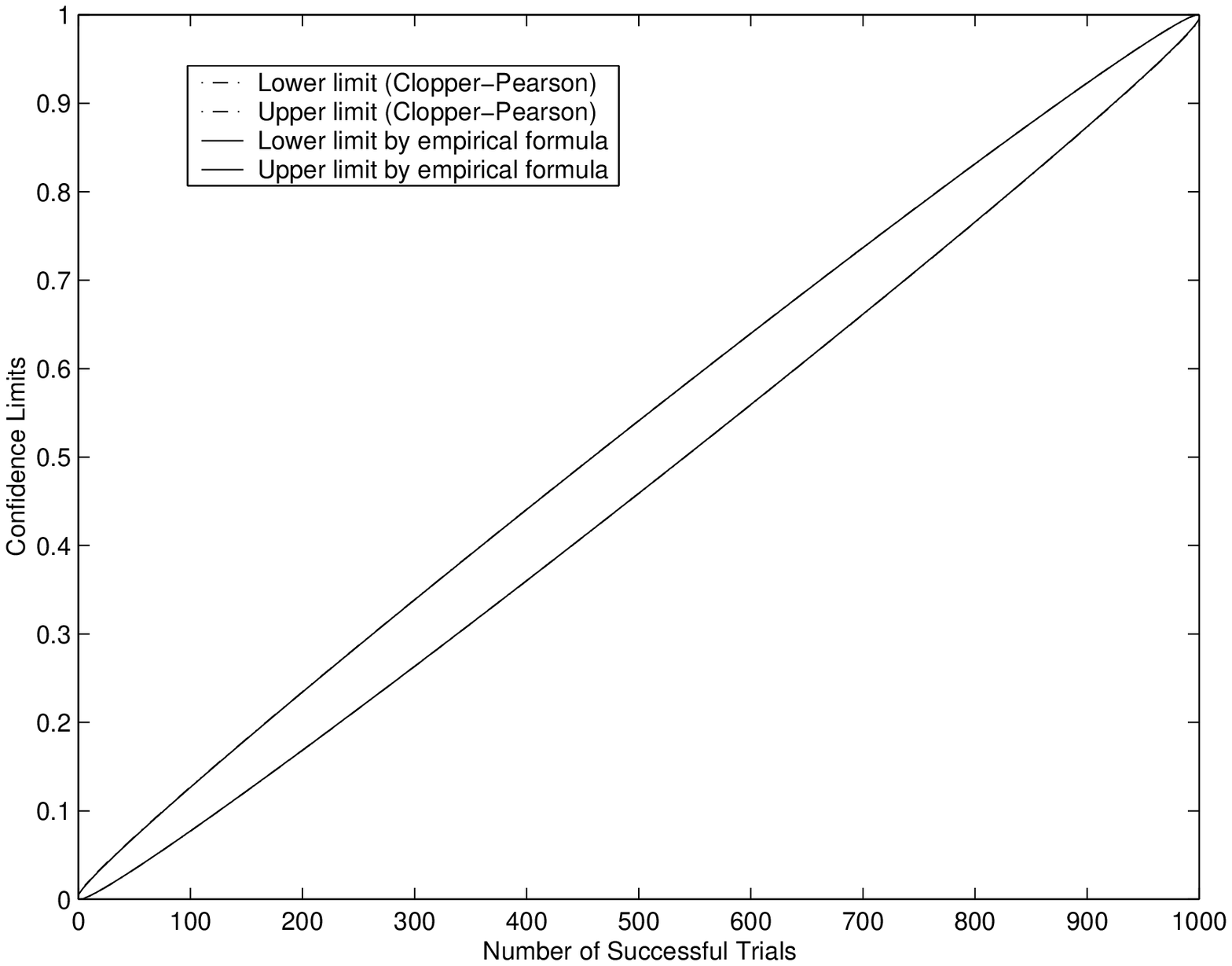, height=3.8in,
width=6.0in }} \caption{ Confidence Interval ($N = 1000, \;
\delta = 0.01$.) } \label{fig_20}
\end{figure}

\begin{figure}[htbp]
\centerline{\psfig{figure= 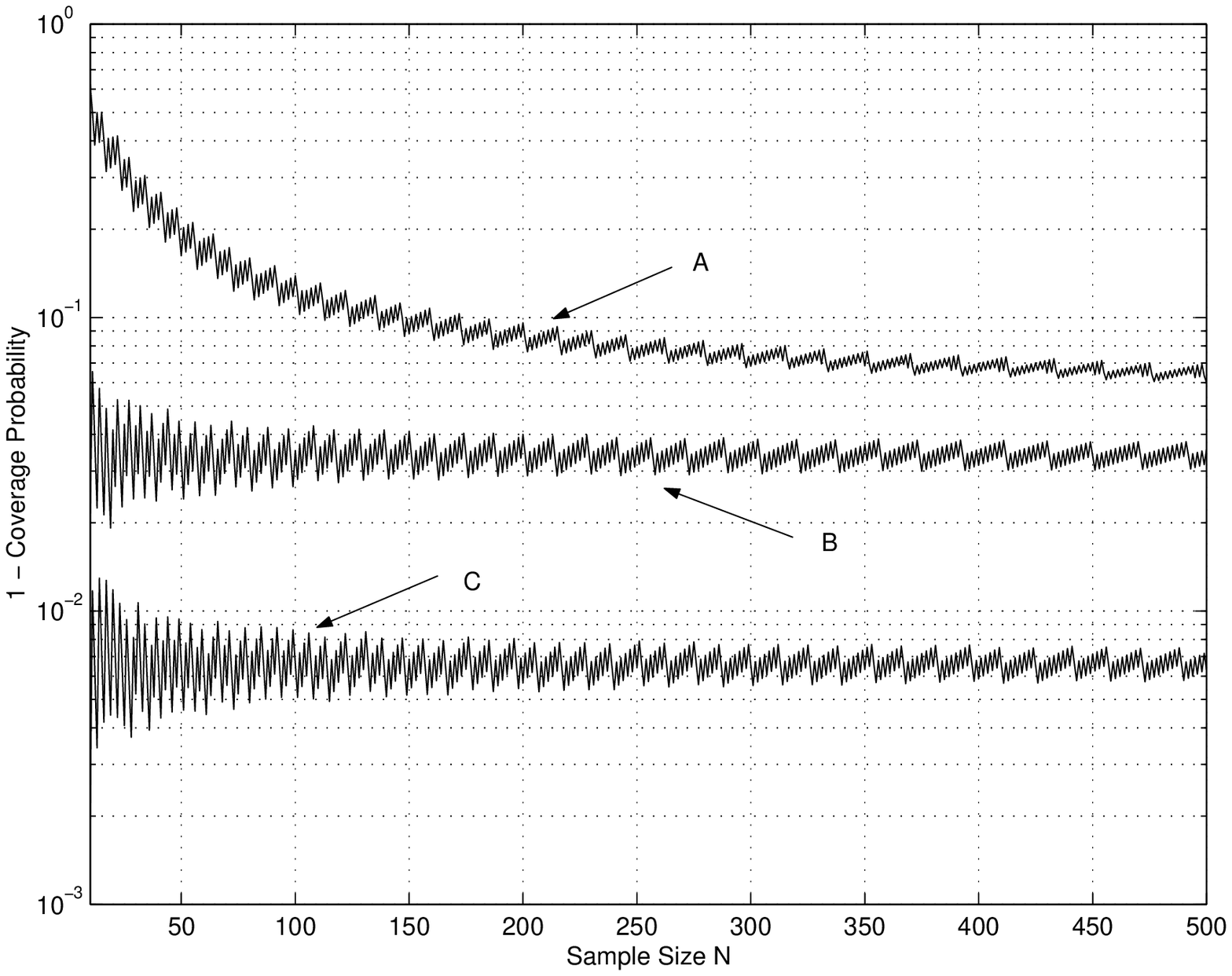, height=3.8in, width=6in
}} \caption{ Error Probability ($\mathbb{P}_X = 0.5, \; \delta =
0.05$. A -- Normal,  B -- Empirical, C -- Rigorous ) }
\label{fig_21}

\bigskip

\bigskip

\centerline{\psfig{figure= 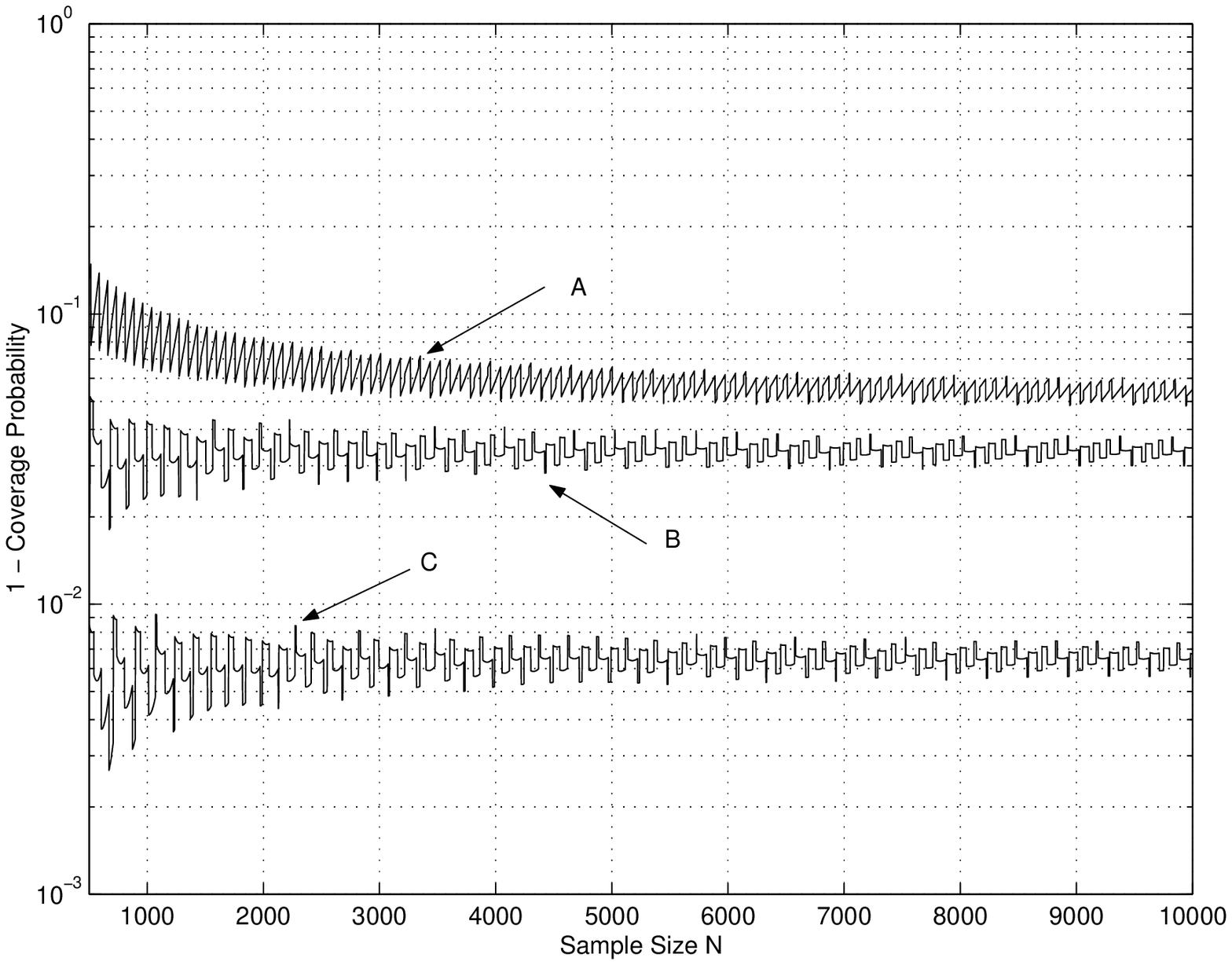, height=3.8in,
width=6in }} \caption{ Error Probability ($\mathbb{P}_X = 0.01, \;
\delta = 0.05$. A -- Normal,  B -- Empirical, C -- Rigorous) }
\label{fig_22}
\end{figure}

\begin{figure}[htbp]
\centerline{\psfig{figure= 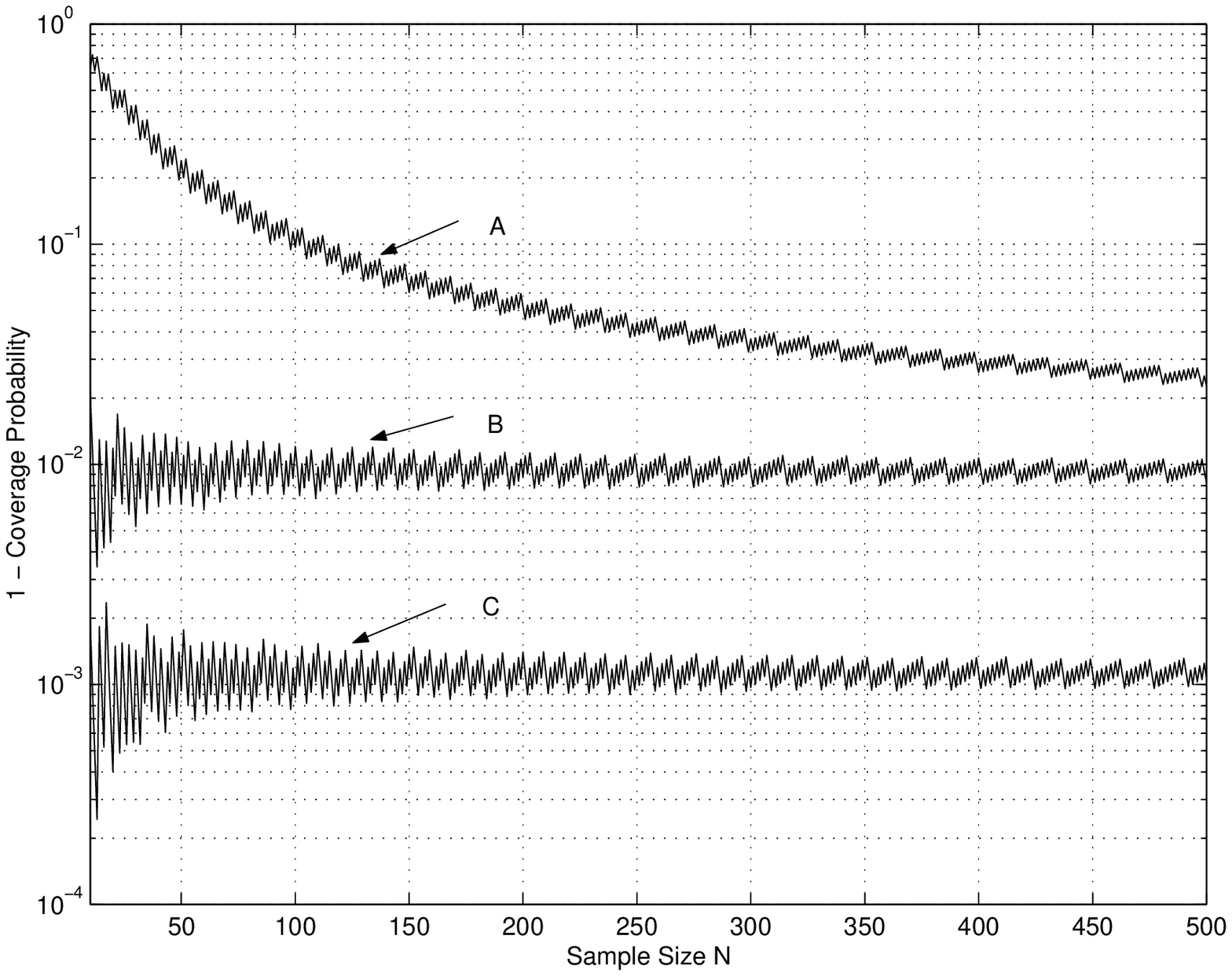, height=3.8in, width=6in
}} \caption{ Error Probability ($\mathbb{P}_X = 0.5, \; \delta =
10^{-2}$. A -- Normal,  B -- Empirical, C -- Rigorous) }
\label{fig_23}

\bigskip

\bigskip

\centerline{\psfig{figure= 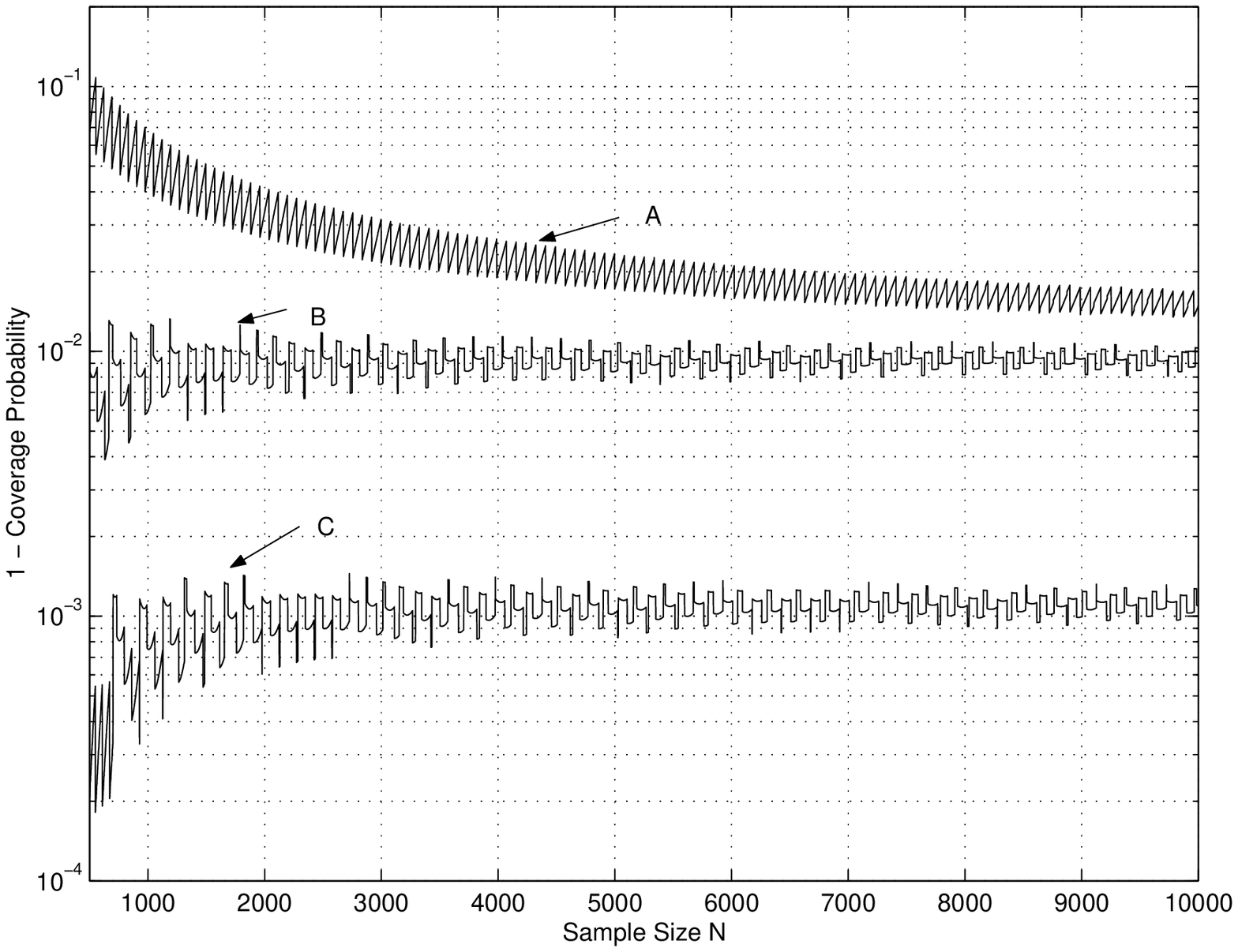, height=3.8in,
width=6in }} \caption{ Error Probability ($\mathbb{P}_X = 0.01, \;
\delta = 10^{-2}$. A -- Normal,  B -- Empirical, C -- Rigorous) }
\label{fig_24}
\end{figure}

\begin{figure}[htbp]
\centerline{\psfig{figure= 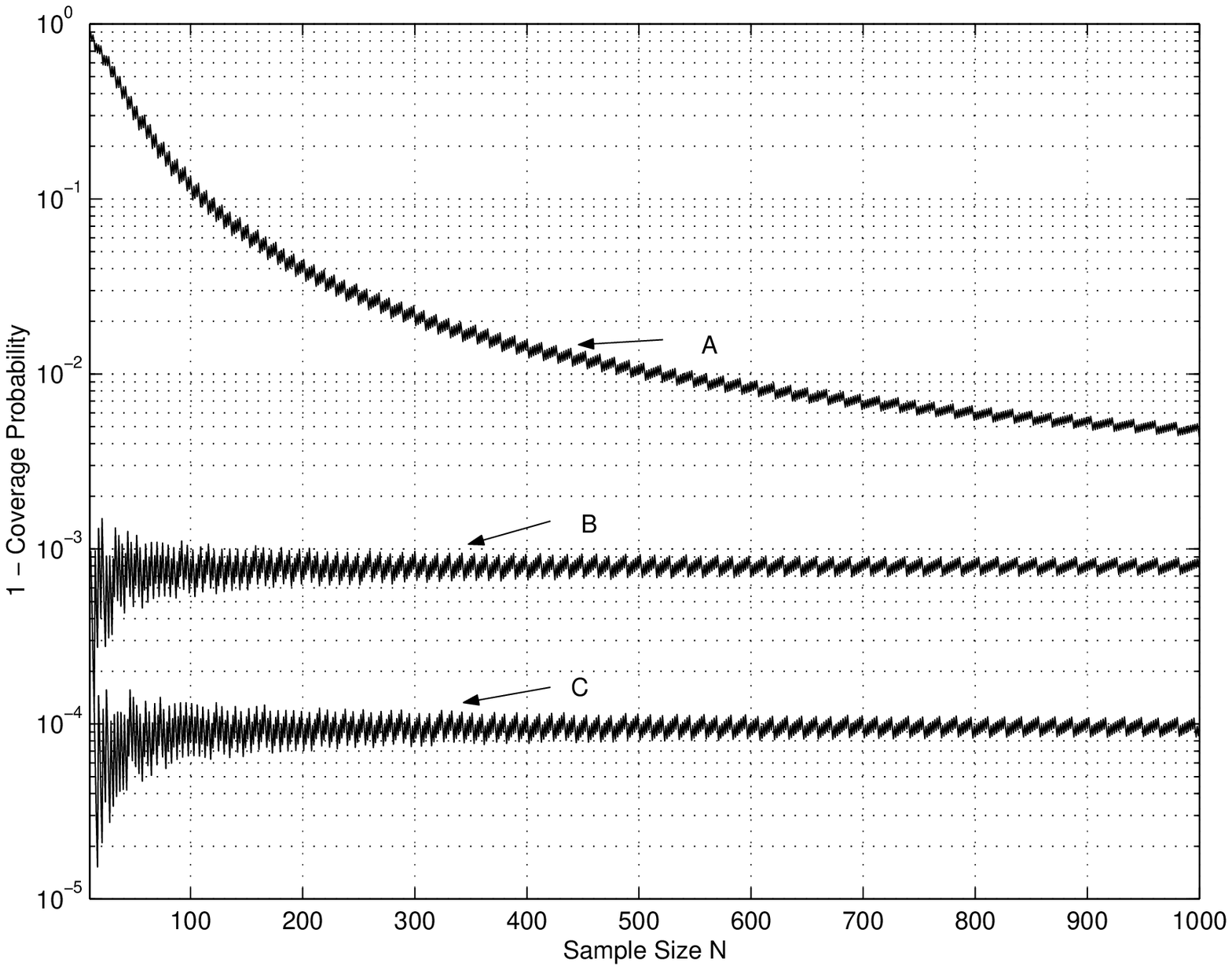, height=3.8in,
width=6in }} \caption{ Error Probability ($\mathbb{P}_X = 0.5, \;
\delta = 10^{-3}$. A -- Normal,  B -- Empirical, C -- Rigorous) }
\label{fig_25}

\bigskip

\bigskip

\centerline{\psfig{figure= 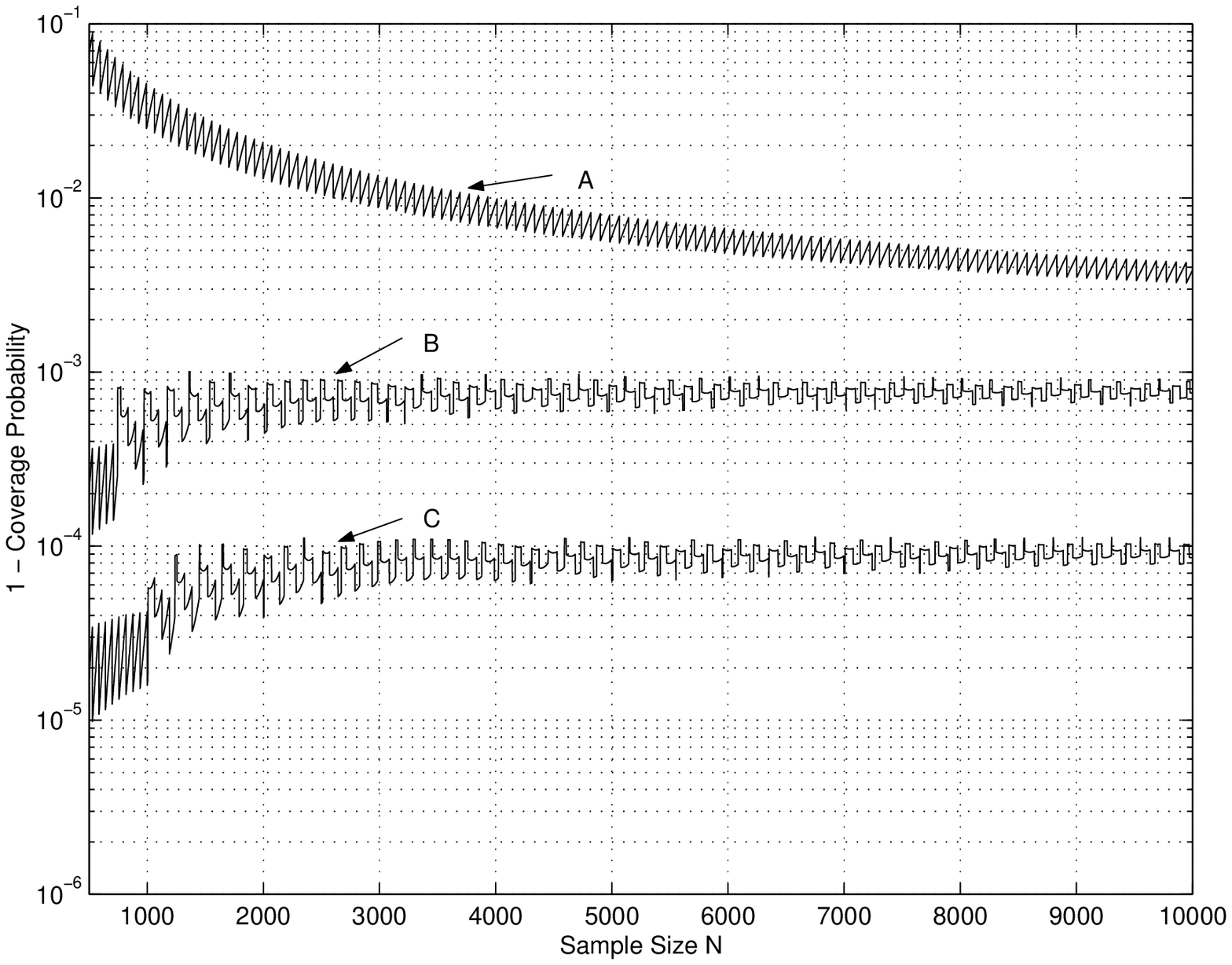, height=3.8in,
width=6in }} \caption{ Error Probability ($\mathbb{P}_X = 10^{-2},
\; \delta = 10^{-3}$. A -- Normal,  B -- Empirical, C --
Rigorous) } \label{fig_26}
\end{figure}

\begin{figure}[htbp]
\centerline{\psfig{figure= 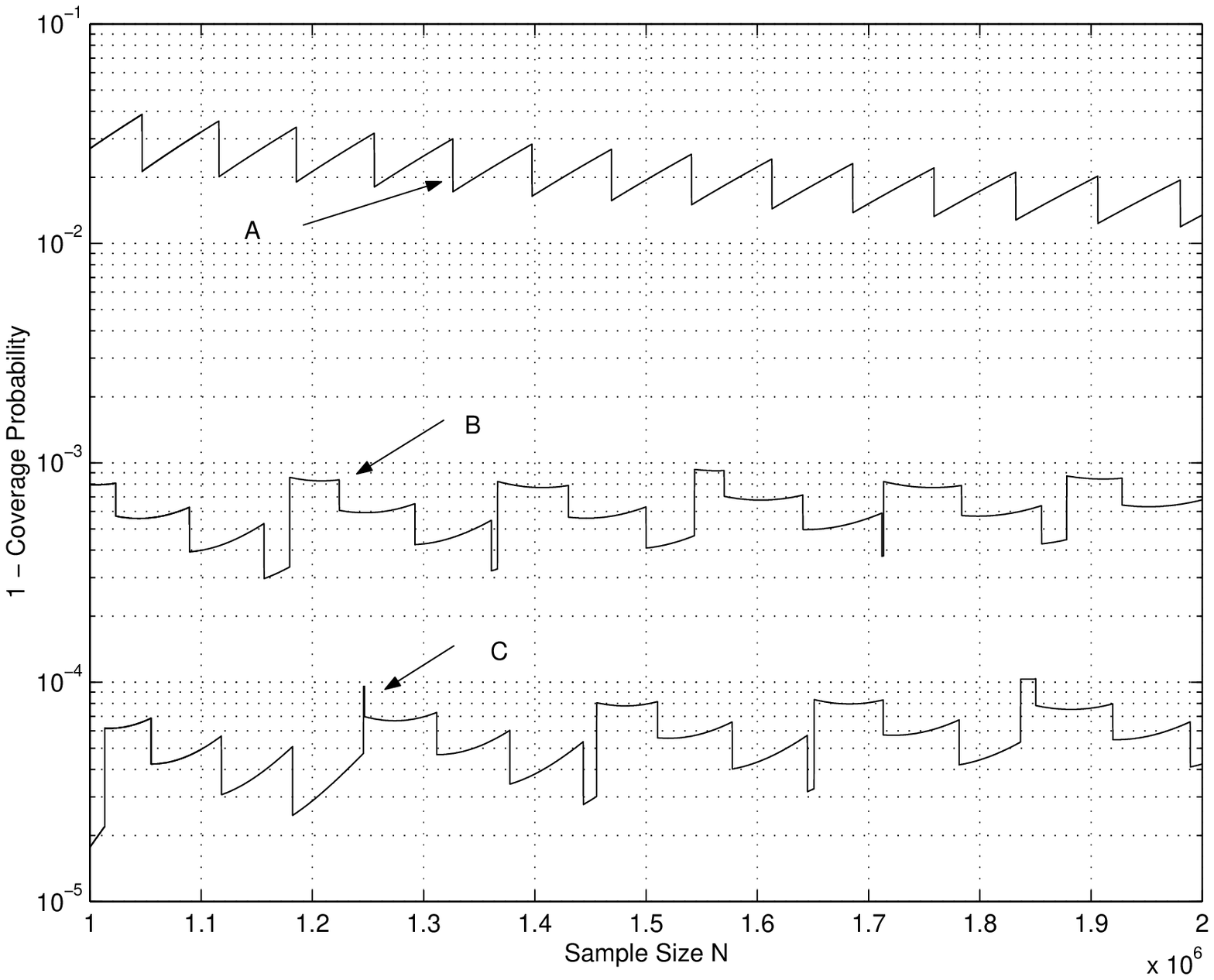, height=3.8in,
width=6in }} \caption{ Error Probability ($\mathbb{P}_X = 10^{-5},
\; \delta = 10^{-3}$. A -- Normal,  B -- Empirical, C --
Rigorous) } \label{fig_27}
\end{figure}

\end{document}